\documentclass[review]{elsarticle}

\usepackage{lineno}
\usepackage{hyperref}
\modulolinenumbers[5]

\RequirePackage{amsmath}
\RequirePackage{amsfonts}
\RequirePackage{bbm}

\newcommand{\Curl}{\ensuremath{\nabla\times}}
\newcommand{\Div}{\ensuremath{\nabla\cdot}}

\newcommand{\permittivity}{\ensuremath{\epsilon}}
\newcommand{\reluctivity} {\ensuremath{\nu}} 

\newcommand{\vecs}[1]{\mathbf{#1}}
 
\newcommand{\Afield}  {\ensuremath{\vecs{A}}} 
\newcommand{\Hcurl}      {\ensuremath{{H}(\mathrm{curl};\Omega)}}

 \journal{Computer Methods in Applied Mechanics and Engineering}

\usepackage[capitalise]{cleveref}
\crefname{equation}{}{}
\usepackage{tikz}
\usepackage{tikz-cd}
\usetikzlibrary{arrows,cd,external}
\usepackage{tikz-3dplot}
\usepackage{pgfplots}
\usepackage{siunitx}
\usepackage{multirow}

\usepackage{subcaption}
\usepackage{nicefrac}

\usepackage[margin=2.5cm]{geometry} \DeclareSymbolFont{Symbols}{OMS}{cmsy}{m}{n}
\SetSymbolFont{Symbols}{bold}{OMS}{cmsy}{b}{n}
\DeclareMathSymbol{\Setminus}{\mathbin}{Symbols}{"6E}
\usepackage{MnSymbol}

\newtheorem{remark}{Remark}
\newtheorem{algor}{Algorithm}

\begin{document}

\begin{frontmatter}

\title{Tree-Cotree Decomposition of Isogeometric Mortared Spaces in H(curl) on Multi-Patch Domains}

\author[epfladdress]{Bernard Kapidani}
\ead{bernard.kapidani@epfl.ch}

\author[cemaddress,cceaddress]{Melina Merkel}
\ead{melina.merkel@tu-darmstadt.de}

\author[cemaddress,cceaddress]{Sebastian Schöps}
\ead{sebastian.schoep@tu-darmstadt.de}

\author[epfladdress,imatiaddress]{Rafael V\'azquez}
\ead{rafael.vazquez@epfl.ch}

\address[epfladdress]{Chair of Numerical Modelling and Simulation, École Polytechnique Fédérale de Lausanne, 1015 Lausanne, Switzerland}
\address[cemaddress]{Computational Electromagnetics Group, Technische Universität Darmstadt, 64289 Darmstadt, Germany}
\address[cceaddress]{Centre for Computational Engineering, Technische Universität Darmstadt, 64289 Darmstadt, Germany}
\address[imatiaddress]{Istituto di Matematica Applicata e Tecnologie Informatiche ``Enrico Magenes'', 27100 Pavia, Italy}

\begin{abstract}
When applying isogeometric analysis to engineering problems, one often deals with multi-patch spline spaces that have incompatible discretisations, e.g. in the case of moving objects. In such cases mortaring has been shown to be advantageous. This contribution discusses the appropriate B-spline spaces needed for the solution of Maxwell's equations in the functions space H(curl) and the corresponding mortar spaces. The main contribution of this paper is to show that in formulations requiring gauging, as in the vector potential formulation of magnetostatic equations, one can remove the discrete kernel subspace from the mortared spaces by the graph-theoretical concept of a tree-cotree decomposition. The tree-cotree decomposition is done based on the control mesh, it works for non-contractible domains, and it can be straightforwardly applied independently of the degree of the B-spline bases. Finally, the simulation workflow is demonstrated using a realistic model of a rotating permanent magnet synchronous machine.
\end{abstract}


\end{frontmatter}

\section{Introduction}

Isogeometric Analysis (IGA) was originally introduced in the context of solid mechanics \cite{Hughes_2005aa}. It is nowadays commonly interpreted as a Finite Element Method (FEM) using spline-based functions for geometry and solution spaces, i.e. generalising the classical polynomial ones. This approach has two advantages, on one hand the geometry can often be exactly described due to the use of the same basis functions (Non-Uniform Rational B-Splines, NURBS) used in computer aided design and on the other hand, the solution is globally more regular, since inter-element smoothness can be controlled by the splines. However, since most computational domains cannot be represented by a single spline-based mapping, multi-patch spaces for geometry and solution have been defined \cite{Buffa_2015aa}. Common multi-patch approaches require consistent spline mappings on both sides of the interface between two patches (i.e. identical knot vectors) and guarantee only continuity across the interface, while smooth multi-patch discretizations need even stronger assumptions, see e.g. \cite{Hughes_2021aa}. The consistency restriction can be mitigated by mortaring, e.g. \cite{Brivadis_2015aa,Buffa_2020aa} at the cost of ensuring the continuity across the interface only weakly.

In the mechanics community nodal element spaces containing square-integrable gradients were in focus. Later, e.g. due to applications in fluid mechanics and electromagnetism, the development of high-order div-conforming Raviart-Thomas and curl-conforming Nédélec-element-like multi-patch spaces \cite{Buffa_2010aa} was pursued. Recently, also the corresponding trace-spaces were investigated and embedded into an exact discrete (De Rham) sequence known from differential geometry \cite{Buffa_2019ac}. Analogously to the continuous case, the discrete sequence for the spline spaces implies that any square-integrable vector field can be split into the sum of a gradient and a curl. 

This paper answers the application-relevant question: how to remove {from the discrete curl-conforming space the kernel of the curl operator, formed by the subspace spanned by gradient fields and harmonic fields.} We will use the graph-theoretical concept of a tree-cotree decomposition. The tree-cotree technique was originally introduced in the finite element community in \cite{Albanese_1988aa} and led to several variants, see the textbook of Bossavit \cite[Section 5.3]{Bossavit_1998aa} and the overview article \cite{Munteanu_2002aa} and the references therein. In the beginning only lowest order discretisations were considered, which simplifies the identification of degrees of freedoms with tree edges. In \cite[Section 6.3]{Zaglmayr_2006aa} it is discussed that the same idea is applicable to the lowest-order subspace of a hierarchical high-order finite-element space. In \cite{Rodriguez_2020aa} the authors investigate a more general tree-cotree decomposition for conventional high-order finite elements on tetrahedra. In this paper we show that the tree for a multi-patch discretisation, which is either glued strongly or weakly by mortaring, can be linked to the tree generated on the control mesh, and this remains valid also for non-contractible domains.

\begin{figure}[t]
  \centering
  \includegraphics[width=.6\textwidth]{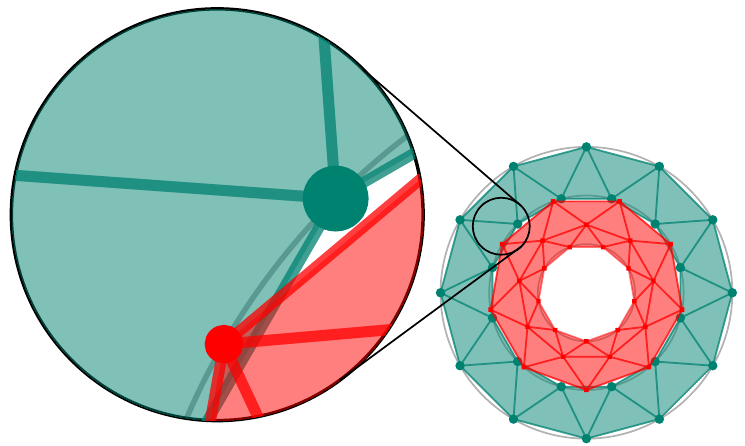}
  \caption{Non-matching rotor/stator interface discretisations when using conventional finite elements.}
  \label{fig:images_mortar_interface}
\end{figure}

{The removal of the discrete kernel of the curl} is a necessary requirement for the well-posedness of many electromagnetic field simulation problems. In particular, the resulting system matrix stemming from discretisation of the curl-curl operator will be singular. This is for example relevant for the magnetostatic curl-curl formulation that is commonly used to simulate electric machines. We demonstrate {the use of the tree-cotree decomposition on} this application using a realistic volumetric-spline model of a rotating permanent magnet synchronous machine. While the formulation and FEM discretisation of such models based on the so-called eddy current problem, e.g. \cite{Biro_1999aa}, is well understood, the application of IGA is recent, \cite{Bontinck_2018ac,Bhat_2018aa,Friedrich_2020aa,Merkel_2021ab}. However, the NURBS-based geometry discretisation promises a robust numerical method since rotor and stator domains share an exact circular interface whose non-matching knot vectors can be overcome by mortaring, see \cref{fig:images_mortar_interface}. 

The paper is structured as follows: \cref{sec:maxwell} introduces Maxwell's equations and the relevant model problems in their strong and weak form. In \cref{sec:iga} isogeometric analysis is summarised, the B-spline spaces and their traces are repeated; finally the mortaring of multi-patch spaces in H(curl) is discussed. \cref{sec:treecotree} constitutes the main contribution, it introduces tree-cotree gauging for the mortared spline spaces. Finally numerical results, i.e. an inf-sup test, eigenvalue and source problems, and an electric machine model, are demonstrated in \cref{sec:numerics} to support the theoretical findings. The paper finishes with conclusions in \cref{sec:conclusions}.

\section{Maxwell Equations for Magnetostatics}\label{sec:maxwell}
Electromagnetic phenomena in the frequency domain are described by Maxwell's equations, i.e.
\begin{align}
 \Curl\vecs{E}&=-j\omega\vecs{B},  
 &\Curl\vecs{H}&=j\omega\vecs{D}+\vecs{J},\label{eq:faradayamperediff}\\
\Div\vecs{D}&=\rho ,
&\Div\vecs{B}&=0,\label{eq:gaussmgaussdiff}
\end{align}
where $\vecs{E}$, $\vecs{D}$, $\vecs{J}$ and $\rho$ are the electric field, displacement, current and charge density, respectively while $\vecs{H}$ and $\vecs{B}$ are the magnetic field and magnetic induction, respectively. Finally $\omega$ is the angular frequency. Maxwell's equations are completed by the material relations 
\begin{align}
 \vecs{D}=\permittivity\vecs{E}, \qquad \vecs{J}=\vecs{J}_{\mathrm{src}}+\sigma\vecs{E}, \qquad
 \vecs{H}=\reluctivity\vecs{B}-\vecs{M},\label{eq:matkap} \end{align}
where $\permittivity$ is the permittivity, $\sigma$ is the electrical conductivity and $\reluctivity$ is the magnetic reluctivity and $\vecs{J}_{\mathrm{src}}$ and $\vecs{M}$ denote the source current density and magnetisation, respectively~\cite{Jackson_1998aa}. All quantities may depend on space and merely for convenience of notation, we disregard field dependencies, e.g. magnetic saturation. 
\subsection{Magnetostatic Source Problem}\label{sec:magnetostatic}
In applications with slowly varying fields one may often assume $\omega=0$. In this case \cref{eq:faradayamperediff,eq:gaussmgaussdiff} decouple into several static cases. We are interested in the magnetostatic setting, e.g. \cite[Chapter 5]{Jackson_1998aa}. This setting is for example particularly relevant for the simulation of permanent magnet synchronous machines, see e.g. \cite{Salon_1995aa}. 
We use the formulation based on the magnetic vector potential $\vecs{A}$ which is introduced such that
\begin{align}
 \vecs{B} = \Curl\vecs{A}. \label{eq:vecpot}
\end{align}

For a bounded domain $\Omega$ with boundary $\partial\Omega$, using \cref{eq:vecpot}, the magnetostatic problem can be derived from Maxwell's equations \cref{eq:faradayamperediff,eq:gaussmgaussdiff} as 
\begin{align}
 \Curl \left( \reluctivity \Curl \vecs{A} \right)&= \vecs{J}_{\mathrm{src}} + \Curl \vecs{M}\qquad &\text{ in } \Omega, \label{eq:magpoteq}\\
            (\Curl \vecs{A}) \times \vecs{n} &= 0 & \text{ on } \Gamma_{N}, \label{eq:neumann} \\
            \vecs{A} \times \vecs{n} &= 0 & \text{ on } \Gamma_{D}, \label{eq:dirichlet}
\end{align}
where $\vecs{n}$ is the outward-pointing normal vector, while $\Gamma_{D}$ and $\Gamma_{N}$ are Dirichlet and Neumann boundaries, respectively ,where $\Gamma_{D} \cup \Gamma_{N} = \partial\Omega$.
The weak formulation of the magnetostatic problem \cref{eq:magpoteq,eq:neumann,eq:dirichlet} can be formulated as:
find $\vecs{A} \in V$ such that
\begin{align}
    \int_{\Omega} \reluctivity (\Curl \vecs{A}) \cdot \left(\Curl \vecs{v}\right) \operatorname{d}V {\;-\int_{\partial\Omega} \reluctivity \left(\Curl \vecs{A} \times \vecs{n}\right)  \cdot \vecs{v} \operatorname{d}s} = \int_{\Omega} \left(\vecs{J}_{\mathrm{src}} + \Curl \vecs{M} \right) \cdot \vecs{v} \operatorname{d}V \label{eq:magnetostatic_weakform}
\end{align}
for all $\vecs{v} \in V$, where $\vecs{v}$ are test functions and $V \subset \Hcurl$ satisfies the boundary condition \eqref{eq:dirichlet}.

If the domain is decomposed into two non-overlapping domains $\overline{\Omega} = \overline{\Omega}_{1} \cup \overline{\Omega}_{2}$ with one common interface $\Gamma_{\mathrm{int}} = \partial \Omega_{1} \cap \partial \Omega_{2}$, the two subdomains can be coupled using mortaring. To do so, for problem \cref{eq:magnetostatic_weakform} we introduce the Lagrange multiplier
\begin{align}
 \boldsymbol{\lambda} = \reluctivity (\Curl \vecs{A}) \times \vecs{n}_{\Gamma},
\end{align}
where $\vecs{n}_{\Gamma}$ is the normal vector on the interface $\Gamma_{\mathrm{int}}$. This results in the new problem:\\
Find $\vecs{A}_{1} \in V_{\vphantom{h,}1}, \vecs{A}_{\vphantom{h,}2} \in V_{\vphantom{h,}2}, \boldsymbol{\lambda} \in M_{\vphantom{h}} $, such that 
\begin{align}
    \int_{\Omega_{1}} \reluctivity (\Curl \vecs{A}_{1}) \cdot (\Curl \vecs{v}_{1}) \operatorname{d}V -\int_{\Gamma_{\mathrm{int}}} \boldsymbol{\lambda}  \cdot \vecs{v}_{1} \operatorname{d}s &= \int_{\Omega_{1}} (\vecs{J}_{\mathrm{src}} + \Curl \vecs{M}) \cdot \vecs{v}_{1} \operatorname{d}V, \label{eq:mortar1}\\
    \int_{\Omega_{2}} \reluctivity (\Curl \vecs{A}_{2}) \cdot (\Curl \vecs{v}_{2}) \operatorname{d}V +\int_{\Gamma_{\mathrm{int}}} \boldsymbol{\lambda}  \cdot \vecs{v}_{2} \operatorname{d}s &= \int_{\Omega_{2}} (\vecs{J}_{\mathrm{src}} + \Curl \vecs{M}) \cdot \vecs{v}_{2} \operatorname{d}V, \label{eq:mortar2}\\
    \int_{\Gamma_{\mathrm{int}}} (\vecs{A}_{1}-\vecs{A}_{2}) \cdot \boldsymbol{\mu} \operatorname{d}s &= 0, \label{eq:mortarcontinuity}
\end{align}
for all $\vecs{v}_{1} \in V_{\vphantom{h,}1}, \vecs{v}_{2} \in V_{\vphantom{h,}2}$ $,  \boldsymbol{\mu} \in M_{\vphantom{h}}$, with the spaces $V_1, V_2 \subset \Hcurl$ satisfying Dirichlet boundary conditions on $\partial \Omega \cap \Gamma_D$, and the space of the Lagrange multiplier $M = \{(\Curl \vecs{v}) \times \vecs{n}_{\Gamma}: \vecs{v} \in V_1\}$.

It is well known, e.g. \cite{Biro_2007aa}, that the weak formulations \eqref{eq:magnetostatic_weakform} and \eqref{eq:mortar1}-\eqref{eq:mortarcontinuity} are not well posed, as the magnetic vector potential $\vecs{A}$ is not unique.
In particular, we can add any gradient to $\vecs{A}$ and, since 
\begin{align}
	\Curl\vecs{A} = \Curl(\vecs{A}+\nabla\xi),
\end{align}
the modified vector potential would still be a valid solution for all $\xi$ with suitable boundary conditions. It is therefore necessary to add a gauging condition to the problem to recover a unique solution.
To this purpose, we will introduce tree-cotree gauging for the discrete isogeometric formulation of those problems.

\subsection{Eigenvalue Problem}
When interested in resonance phenomena on a domain $\Omega$, one neglects magnetisation and current densities and rearranges \eqref{eq:faradayamperediff}-\eqref{eq:gaussmgaussdiff} such that one obtains Maxwell's eigenvalue problem in terms of the electric field strength, i.e.
\begin{align}
	\Curl\left(\reluctivity\Curl{\vecs{E}}\right) &= \omega^{2}\permittivity{\vecs{E}}
	\qquad &\text{ in } \Omega,
	\label{eq:eigvalprob}
\end{align}
with unknown eigenmodes $\vecs{E}$ and eigenvalues $\omega$. We will use this eigenvalue problem for the numerical study of the method, for which we will always consider the Dirichlet condition $\vecs{E} \times \vecs{n}=0$ on $\partial\Omega$. We here omit the corresponding weak and mortaring formulation of \eqref{eq:eigvalprob} since they are well known and similar to the magnetostatic case.

\section{Isogeometric Analysis and Mortar Method}\label{sec:iga}
In the first part of this section we define spline spaces for the discretisation of the Maxwell equations, mainly to introduce notation, and define a suitable space for the Lagrange multiplier for mortaring in the single-patch case, following previous work \cite{Buffa_2020aa}. In the second part we introduce a modification for the multi-patch case, which will prove to be crucial for applying the tree-cotree gauging.
\subsection{Isogeometric Mortar Formulation on Single Patch Subdomains}
\subsubsection{B-splines and NURBS}
A set of $n$ univariate B-splines of degree $p$ is constructed from an ordered knot vector $\Xi = \{\xi_1, \ldots, \xi_{n+p+1}\}$, with $0 \le \xi_i \le \xi_{i+1} \le 1$ for every $i$. We assume that the knot vector is open, which means that the first and last knot are repeated exactly $p+1$ times. We denote by $\hat B_i^p$ the $i$-th basis function of degree $p$ in the reference domain $(0,1)$, and by $S_p(\Xi)$ the space spanned by these B-splines.

Multivariate spaces in the reference domain $(0,1)^d$, usually with $d=2,3$, are defined by tensor-product of univariate B-splines. Introducing a degree $p_j$ and a knot vector $\Xi_j$ in the $j$th coordinate, for $j=1, \ldots, d$, the multivariate basis functions are defined as
\[
\hat B_{\mathbf{i}}^{\mathbf{p}} (\boldsymbol{\xi}) = \hat B_{i_1}^{p_1}(\xi_1) \ldots \hat B_{i_d}^{p_d}(\xi_d),
\]
with vector degree $\mathbf{p} = (p_1, \ldots, p_d)$, and multivariate index $\mathbf{i} = (i_1, \ldots, i_d)$. They span the multivariate spline space $S_{p_1,\ldots,p_d}(\Xi_1,\ldots,\Xi_d) = S_{p_1}(\Xi_1) \otimes \ldots \otimes S_{p_d}(\Xi_d)$. NURBS basis functions are defined from B-splines by associating a weight to each B-spline function, see \cite{Hughes_2005aa} for details.

\begin{figure}
    \centering
\includegraphics[width=0.4\textwidth,trim=2cm 1cm 2cm 0cm,clip]{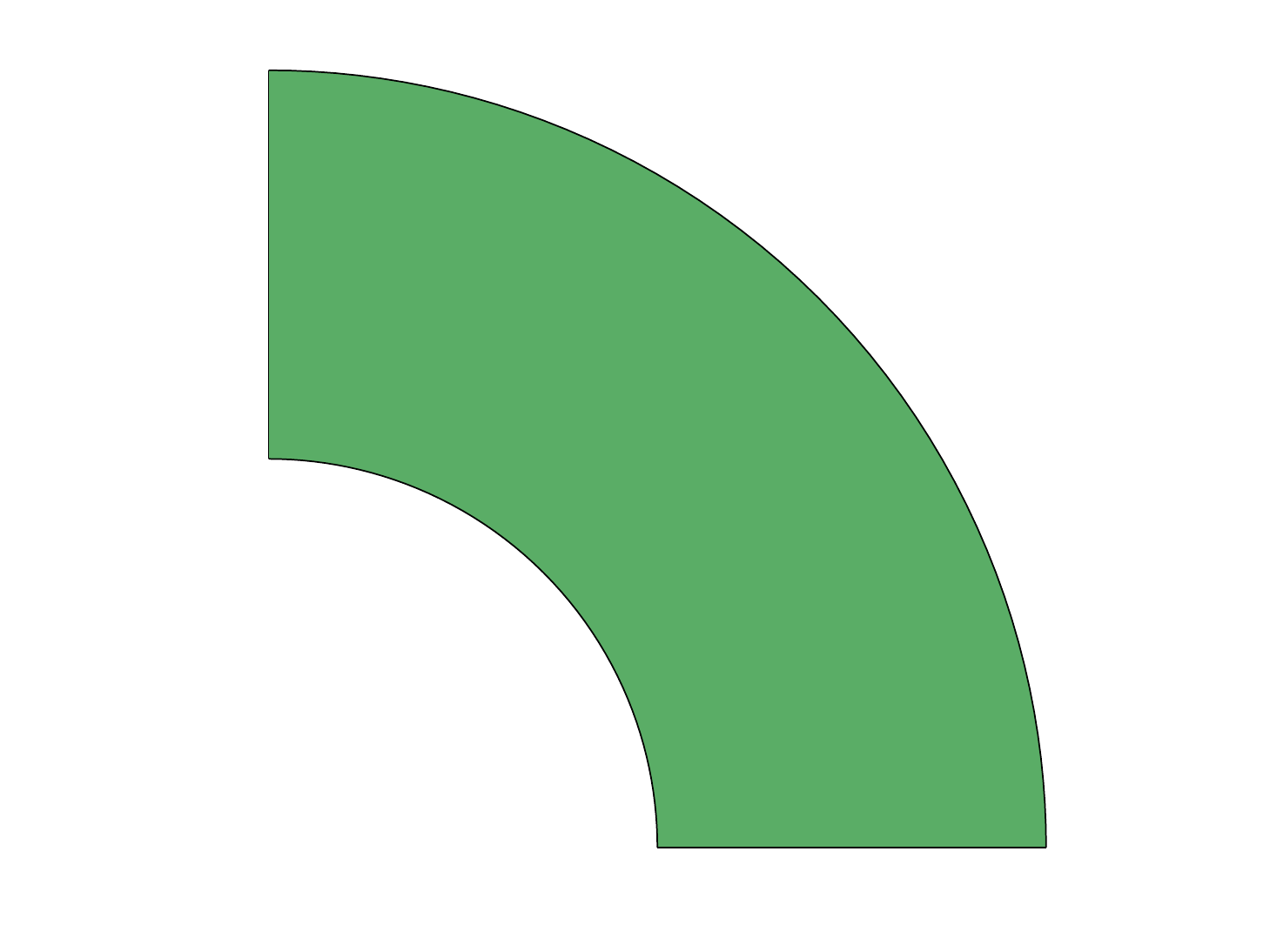}
        \includegraphics[width=0.4\textwidth,trim=2cm 1cm 2cm 0cm,clip]{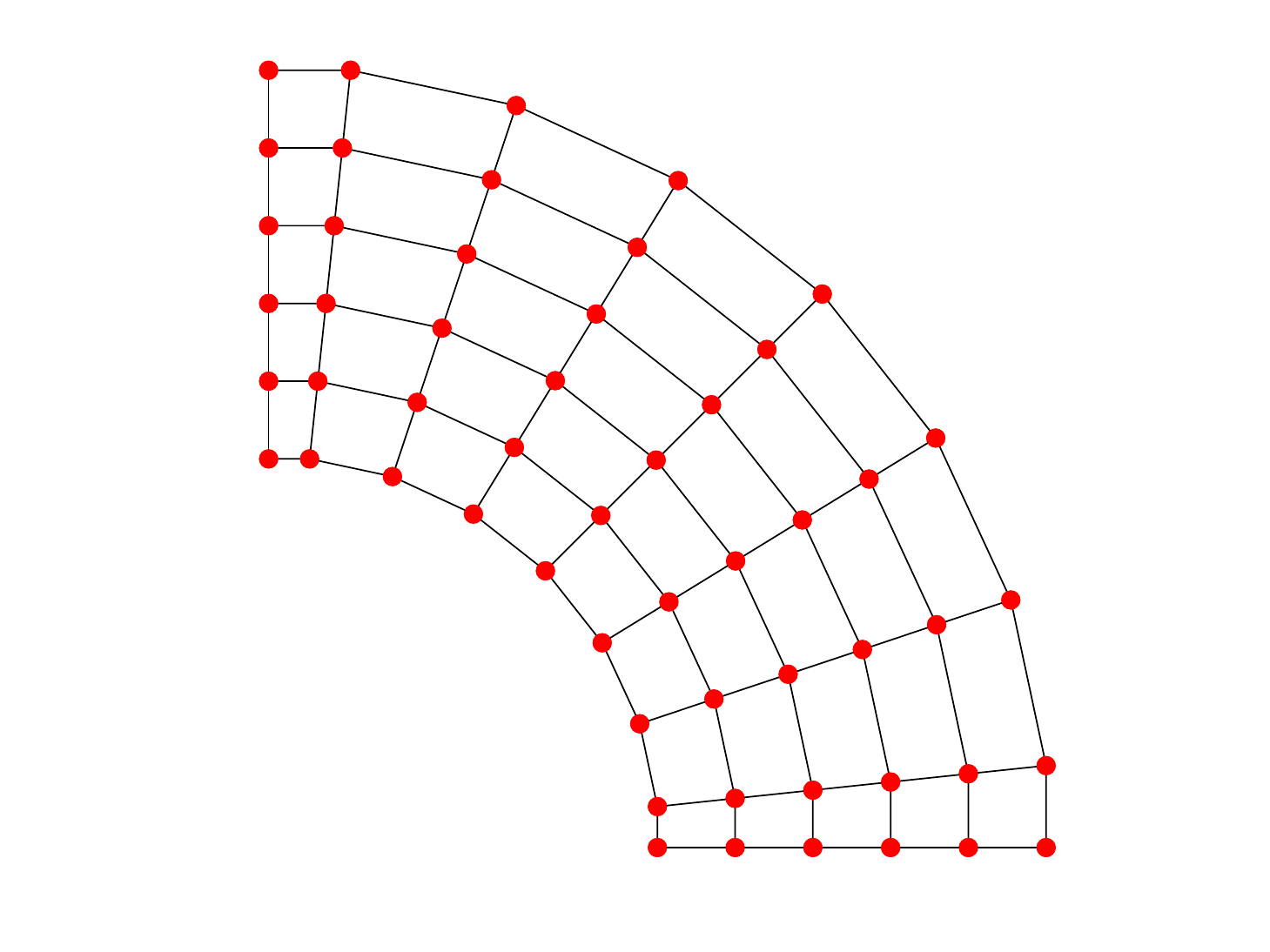}
\caption{A quarter ring geometry, defined as the set $[1/2,1] \times [0,\pi/2]$ in the $(r,\varphi)$ polar coordinates, and the control mesh of its NURBS representation.} \label{fig:ctrl_mesh_example}
\end{figure}

A spline geometry is built as a map from the reference domain $(0,1)^d$ to the physical domain by associating a control point $\mathbf{P}_{\mathbf{i}} \in \mathbb{R}^k$ to each basis function, namely
\[
\mathbf{F}(\boldsymbol{\xi}) = \sum_{\mathbf{i}} {\mathbf{P}}_{\mathbf{i}} \hat B_{\mathbf{i}}^{\mathbf{p}} (\boldsymbol{\xi}), \quad \boldsymbol{\xi} \in (0,1)^d,
\]
which defines a map from $\hat \Omega = (0,1)^d$ into $\Omega \subset \mathbb{R}^k$. The set of control points determines a structured control mesh, see \cref{fig:ctrl_mesh_example}. NURBS parametrisations are defined analogously, replacing B-spline functions by NURBS basis functions.

\subsubsection{Spline Spaces for the Maxwell Equations}
From now on we assume that the domain $\Omega \subset \mathbb{R}^3$ is defined through a trivariate NURBS mapping, and that it is a bi-Lipschitz homeomorphism, to prevent singularities, see \cite[Section~3]{Beirao-da-Veiga_2014aa}. We also assume for simplicity that the degrees are equal in all directions, that is, $p_j = p$ for every $j$.

Starting from the knot vector $\Xi_j$, we define $\Xi'_j$ as the modified knot vector removing the first and last knots, and also the univariate spaces $S_{p-1}(\Xi_j')$ of degree $p-1$. For these spaces we replace the standard basis by Curry-Schoenberg B-splines $\hat D_{i,p-1} = \frac{p}{\xi_{i+p} - \xi_{i}} \hat B_{i,p-1}$, which yields the following expression for the derivative of degree $p$ splines
\begin{equation} \label{eq:derivative}
\hat B'_{i,p} = \hat D_{i,p-1} - \hat D_{i+1,p-1}.
\end{equation}
Based on these univariate spaces, basis functions, and on the parametrisation $\mathbf{F}$, following \cite[Section~5]{Beirao-da-Veiga_2014aa} we define a spline De Rham complex
\begin{equation}
\begin{tikzcd}
S_p^0(\Omega) \arrow[r,"\nabla"] & S_p^1(\Omega) \arrow[r,"\nabla \times"] & S_p^2(\Omega) \arrow[r,"\nabla \cdot"] & S_p^3(\Omega),
\end{tikzcd}\label{eq:spline_derham}
\end{equation}
where in particular the curl-conforming space is given by
\begin{align*}
& S_p^1(\Omega) = \{ \mathbf{u} \in \Hcurl : \mathbf{u} \circ \mathbf{F} = (D \mathbf{F}^\top)^{-1} \hat {\mathbf{u}}, \, \hat {\mathbf{u}} \in S_p^1(\hat \Omega) \}, \\
& S_p^1(\hat \Omega) = S_{p-1,p,p}(\Xi'_1,\Xi_2,\Xi_3) \times S_{p,p-1,p}(\Xi_1,\Xi'_2,\Xi_3) \times S_{p,p,p-1}(\Xi_1,\Xi_2,\Xi'_3),
\end{align*}
and $D\mathbf{F}$ is the Jacobian of the parametrisation.

To impose Dirichlet boundary conditions, let $\Sigma \subseteq \partial \Omega$, and for simplicity we assume that $\Sigma$ is the image through $\mathbf{F}$ of some boundary sides of the reference domain. We have an analogous diagram for spaces with vanishing boundary conditions on $\Sigma$, given by
\begin{equation*}
\begin{tikzcd}
S_p^0(\Omega;\Sigma) \arrow[r,"\nabla"] & S_p^1(\Omega;\Sigma) \arrow[r,"\nabla \times"] & S_p^2(\Omega;\Sigma) \arrow[r,"\nabla \cdot"] & S_p^3(\Omega;\Sigma).
\end{tikzcd}
\end{equation*}
In particular, the curl-conforming space is 
\begin{align*}
& S_p^1(\Omega;\Sigma) = \{ \mathbf{u} \in S_p^1(\Omega) : (\mathbf{n} \times \mathbf{u}) \times \mathbf{n} = \mathbf{0} \text{ on } \Sigma \}.
\end{align*}

It is important to note that the basis functions of $S_p^0(\Omega)$ are associated to the control points of the control mesh, and analogously, the basis functions of $S_p^1(\Omega)$ are associated to the edges of the control mesh, analogously to the classical Nédélec basis. Similarly, basis functions for the spaces with vanishing boundary conditions are obtained removing the degrees of freedom associated to the control points or edges corresponding to $\Sigma$. Moreover, and thanks to \eqref{eq:derivative}, the derivatives of the basis functions of $S_p^0(\Omega)$ are written in terms of the basis functions of $S_p^1(\Omega)$, with coefficients given by the vertex-edge incidence matrix of the control mesh. All these properties are summarised by the fact that there exists a set of isomorphisms relating the spline spaces to low-order finite element spaces in the control mesh, and these isomorphisms commute with the differential operators, see \cite{BSV13} and \cite[Chapter 5]{Beirao-da-Veiga_2014aa} for more details.

\subsubsection{Trace Spaces} \label{sec:traces}
Following \cite{Buffa_2020aa}, we let $\Gamma \subset \partial \Omega$, and we assume that $\Gamma$ is the image through $\mathbf{F}$ of one boundary side of the reference domain. Therefore, there exists $\mathbf{F}_\Gamma: \hat\Gamma = (0,1)^2 \rightarrow \Gamma \subset \mathbb{R}^3$. Then, on the surface $\Gamma$ we can define two different sequences of spline spaces
\begin{equation*}
\begin{tikzcd}
S_p^0(\Gamma) \arrow[r,"\nabla"] & S^1_p(\Gamma) \arrow[r,"\nabla \times"] & S^2_p(\Gamma), \\
S_p^0(\Gamma) \arrow[r,"\nabla^\perp"] & S^{1^*}_p(\Gamma) \arrow[r,"\nabla \cdot"] & S^2_p(\Gamma),
\end{tikzcd}
\end{equation*}
where in particular $S^1_p(\Gamma)$ and $S^{1^*}_p(\Gamma)$ are respectively curl-conforming and div-conforming spaces, defined from suitable push-forwards from the spaces in the parametric domain $\hat \Gamma = (0,1)^2$ given by 
\begin{align*}
& S_p^1(\hat \Gamma) = S_{p-1,p}(\Xi'_k,\Xi_l) \times S_{p,p-1}(\Xi_k,\Xi'_l), \\
& S_p^{1^*}(\hat \Gamma) = S_{p,p-1}(\Xi_k,\Xi'_l) \times S_{p-1,p}(\Xi'_k,\Xi_l),
\end{align*}
with the indices of the knot vectors $1 \le k < l \le 3$ depending on the position of the boundary side. These spaces are also the image of suitable trace operators applied to the spaces $S_p^i(\Omega)$. In fact, if we define the trace operators
\begin{align*}
\gamma^0_\Gamma(u) =  u|_\Gamma, \qquad
{\gamma}^1_\Gamma (\mathbf{u}) = (\mathbf{n} \times \mathbf{u}) \times \mathbf{n}, \qquad
{\gamma}^{1^*}_{\Gamma} (\mathbf{u}) = \mathbf{u} \times \mathbf{n}, \qquad
\gamma^2_\Gamma (\mathbf{u}) = \mathbf{u} \cdot \mathbf{n},
\end{align*}
it holds that
\begin{align*}
\gamma^0_\Gamma(S^0_p(\Omega)) = S^0_p(\Gamma), \qquad 
\gamma^1_\Gamma(S^1_p(\Omega)) = S^1_p(\Gamma), \qquad 
\gamma^{1^*}_\Gamma(S^1_p(\Omega)) = S^{1^*}_p(\Gamma), \qquad 
\gamma^2_\Gamma(S^2_p(\Omega)) = S^2_p(\Gamma).
\end{align*}
We can also define the spaces with vanishing boundary conditions, and in particular functions in $S_p^0(\Gamma; \partial \Gamma)$ vanish on the boundary, for $S_p^1(\Gamma;\partial \Gamma)$ the tangential trace vanishes, and for $S_p^{1^*}(\Gamma;\partial \Gamma)$ the normal trace vanishes on the boundary of $\Gamma$.

Finally, assuming that $\Lambda$ is a one-dimensional side of $\Gamma$, and therefore parameterised as $\mathbf{F}_\Lambda: (0,1) \rightarrow \Lambda$, we will also make use of the one-dimensional sequence
\begin{equation*}
\begin{tikzcd}
S_p^0(\Lambda) \arrow[r,"\nabla"] & S^1_p(\Lambda),
\end{tikzcd}
\end{equation*}
where the spaces are image through suitable push-forward of the ones defined in the unit interval, $S_p(\Xi)$ and $S_{p-1}(\Xi')$, respectively. Moreover, they can be obtained as trace spaces of $S_p^0(\Gamma)$ and $S_p^1(\Gamma)$. 

\subsubsection{Mortar Formulation}
Let us assume that the two subdomains in formulation \eqref{eq:mortar1}--\eqref{eq:mortarcontinuity}, $\Omega_1$ and $\Omega_2$, are defined by a NURBS parametrisation as described in previous sections, that the interface $\Gamma_{\mathrm{int}}$ corresponds to a boundary side of each subdomain, and that $\Gamma_D = \partial \Omega$. Let us denote by $\Sigma_k = \partial \Omega_k \Setminus \Gamma_{\mathrm{int}}$ for $k=1,2$. Then, the discrete spaces for the approximation of the curl-conforming field are chosen as the adequate spline spaces, that is, $V_{k,h} = S_p^1(\Omega_k;\Sigma_k)$ for $k=1,2$. Note that the two associated meshes are not necessarily conforming on the interface $\Gamma_{\mathrm{int}}$.

The discrete space for the Lagrange multiplier in the mortar method follows the construction introduced for the first time in \cite{Buffa_2020aa}. We build the multiplier space on the $\Omega_1$ subdomain, and as it will constrain the functions in $\Omega_1$, we say that $\Omega_2$ is independent while $\Omega_1$ is the dependent subdomain\footnote{In this paper we prefer the nomenclature of independent/dependent subdomains, however this is also often referred to as master/slave in the literature.}. The construction is simply based on choosing the div-conforming space $M_h = S_{p-1}^{1^*}(\Gamma_{\mathrm{int}})$, i.e. a spline space with one degree lower than the one for the curl-conforming field, and associated to the same mesh (or knot vectors) of the space in the dependent subdomain. It was proved in \cite{Buffa_2020aa} that the space has the same dimension as the trace space $\gamma^1_{\Gamma_{\mathrm{int}}}\left(S_p^1(\Omega_1;\Sigma_1)\right)$, and that it is stable, in the sense that the following inf-sup stability condition is fulfilled
\[
\sup_{\mathbf{u} \in \gamma^1_{\Gamma_{\mathrm{int}}}\left(S_p^1(\Omega_1;\Sigma_1)\right)} \frac{\int_{\Gamma_{\mathrm{int}}} \mathbf{u} \cdot \mathbf{v}\operatorname{d}s}{\parallel \mathbf{u} \parallel_{-1/2,\mathbf{curl}}} \ge \beta \parallel \mathbf{v} \parallel_{-1/2,\mathrm{div}} \quad \text{ for all } \mathbf{v} \in S_{p-1}^{1^*}(\Gamma_{\mathrm{int}}).\label{eq:inf_sup}
\]

In fact, the same properties (equal dimension and inf-sup stability) are also valid for the trace space $\gamma^0_{\Gamma_{\mathrm{int}}} (S_p^0(\Omega_1;\Sigma_1))$ with the multiplier space $S_{p-1}^2(\Gamma_{\mathrm{int}})$, and for $\gamma^2_{\Gamma_{\mathrm{int}}} (S_p^2(\Omega_1;\Sigma_1))$ with $S_{p-1}^0(\Gamma_{\mathrm{int}})$.

\subsection{Isogeometric Mortar Formulation on Multi-Patch Subdomains} \label{sec:mortar-multipatch}
We will now detail the case in which both the independent and the dependent domain are formed by multiple patches, as this will require a modification of the space for the Lagrange multiplier. We will assume conforming patches on each subdomain, in such a way that mortaring only needs to be applied between $\Omega_1$ and $\Omega_2$, but not within them.

\subsubsection{Mortar Space on Multi-Patch Domains: the Ungauged Case} \label{sec:multipatch-ungauged}
Let $\Omega_1 = \cup_{i=1}^{N_1} \Omega^i_1$ and $\Omega_2 = \cup_{i=1}^{N_2} \Omega^i_2$, each $\Omega^i_k$ being the image of a parametrisation $\mathbf{F}^i_k$. We assume that they all have empty intersections, and denoting the interfaces by $\Gamma^{ij}_k = \partial \Omega_k^i \cap \partial \Omega_k^j$, we assume that they are conforming in the sense that $\Gamma_k^{ij}$, for $k=1,2$ and $i,j = 1, \ldots, N_k$ is either empty or the image of a full side of the parametric domain both for $\mathbf{F}_k^i$ and $\mathbf{F}_k^j$, and the knot vectors and the parametrisation from both sides match (see \cite[Section 3.2]{Beirao-da-Veiga_2014aa} for details). To introduce boundary conditions, we will also denote by $\Sigma^i_k = \partial \Omega_k^i \cap \partial \Omega$.

Leveraging the single patch definitions, we define the discrete curl-conforming space on each multi-patch domain, with vanishing tangential traces on the boundary of $\Omega$, as
\begin{equation*}
V_{k,h} = \{ \mathbf{u}_h \in H(\mathrm{curl};\Omega_k) : \mathbf{u}_h|_{\Omega_k^i} \in S_p^1(\Omega_k^i; \Sigma_k^i) \text{ for } i = 1, \ldots, N_k \},
\end{equation*}
where we note that tangential traces inside $\Omega_k$ are glued together strongly, by identifying the coincident control edges of adjacent patches. Note also that the tangential traces do not vanish on $\Gamma_{\mathrm{int}}$.

The spaces $V_{k,h}$ for $k=1,2$ are used for the discretisation of the curl-conforming field, but we still need to introduce the space $M_h$ for the Lagrange multipler. Assuming for simplicity that each subdomain $\Omega_1^i$ contains only one face on $\Gamma_{\mathrm{int}}$, that we denote by $\Gamma^i$, we define this space as
\[
M_h = \{ \boldsymbol{\mu}_h \in (L^2(\Gamma_{\mathrm{int}}))^3 : \boldsymbol{\mu}_h \in S_{p-1}^{1^*}(\Gamma^i) \text{ for } i = 1, \ldots, N_1 \},
\]
where we note that, unlike the previous spaces, the functions of $M_h$ are not glued together along the interfaces. This space was introduced for mortaring in \cite{Buffa_2020aa} and used for the solution of Maxwell eigenvalue problem, giving a stable and spurious-free discretisation. However, as we will see in the following, it is not a good choice when gauging is needed, because the discrete kernel of the curl operator does not coincide with the gradients of a suitable discrete space.

\subsubsection{Mortar Space on Multi-Patch Domains: Modification for the Gauged Case} \label{sec:mortar-gauged}
To analyse the discrete kernel we consider a simple but analogous problem. We restrict ourselves to the dependent domain $\Omega_1$ and replace the mortar gluing by a Dirichlet problem, with tangential boundary conditions applied through the Lagrange multiplier on $\Gamma_{\mathrm{int}}$, and strongly elsewhere. Note that $\Omega_1$ remains a multi-patch domain. Using the Lagrange multipliers, our discrete space is given by
\[
X^1_{h} = \{ \mathbf{u}_h \in V_{1,h} : \left( (\mathbf{n} \times \mathbf{u}_h) \times \mathbf{n}, \boldsymbol{\mu}_h \right)_{\Gamma_{\mathrm{int}}} = 0 \; \; \forall \boldsymbol{\mu}_h \in M_h \}.
\]
We also define the discrete kernel as the subspace
\[
K_{h} = \{ \mathbf{u}_h \in X^1_h : (\nabla \times \mathbf{u}_h, \nabla \times \mathbf{v}_h)_{\Omega_1} = 0 \; \forall \mathbf{v}_h \in X^1_h \}.
\]
To analyse whether $K_{h}$ contains all the necessary gradients, we introduce the spaces of $H^1$-conforming functions without and with vanishing boundary conditions on $\Gamma_{\mathrm{int}}$:
\begin{align*}
& X^0_{h} = \{ \phi_h \in H^1(\Omega_1) : \phi_h|_{\Omega_1^i} \in S_p^0(\Omega_1^i;\Sigma_1^i) \text{ for } i = 1, \ldots, N_1 \}, \\
& X^0_{h,0} = \{ \phi_h \in H_0^1(\Omega_1) : \phi_h|_{\Omega_1^i} \in S_p^0(\Omega_1^i;\Sigma_1^i) \text{ for } i = 1, \ldots, N_1 \} = X^0_{h} \cap H^1_0(\Omega_1).
\end{align*}
To have a correct discrete kernel, we should have $\nabla X^0_{h,0} = K_{h}$. While it is easy to prove that $\nabla X^0_{h,0} \subset K_{h}$, in general the equality is not true. In fact, if we denote by $Z_{\Gamma_{\mathrm{int}}}$ the set of vertices between patches internal to $\Gamma_{\mathrm{int}}$, we conjecture that the dimension of the discrete kernel is
\[
\dim X^0_{h,0} \le \dim K_h = \dim X^0_{h,0} + \# Z_{\Gamma_{\mathrm{int}}},\]
which we here show with a numerical example. A more rigorous discussion about this conjecture is given in \ref{sec:appendix-kernel}.

\begin{figure}
\centering
\includegraphics[width=0.3\textwidth,trim=13cm 0cm 11cm 0cm, clip]{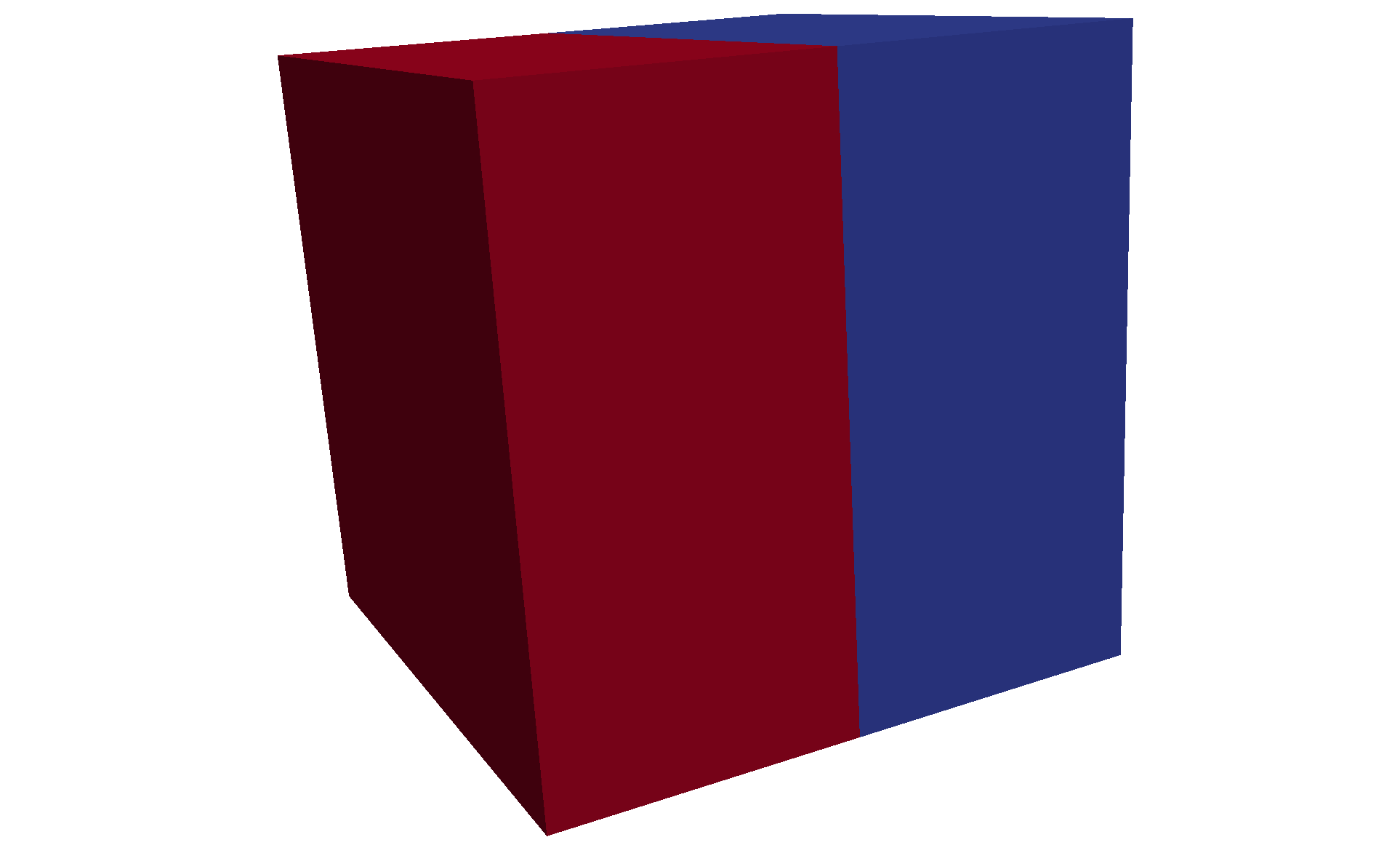}
\includegraphics[width=0.3\textwidth,trim=13cm 0cm 11cm 0cm, clip]{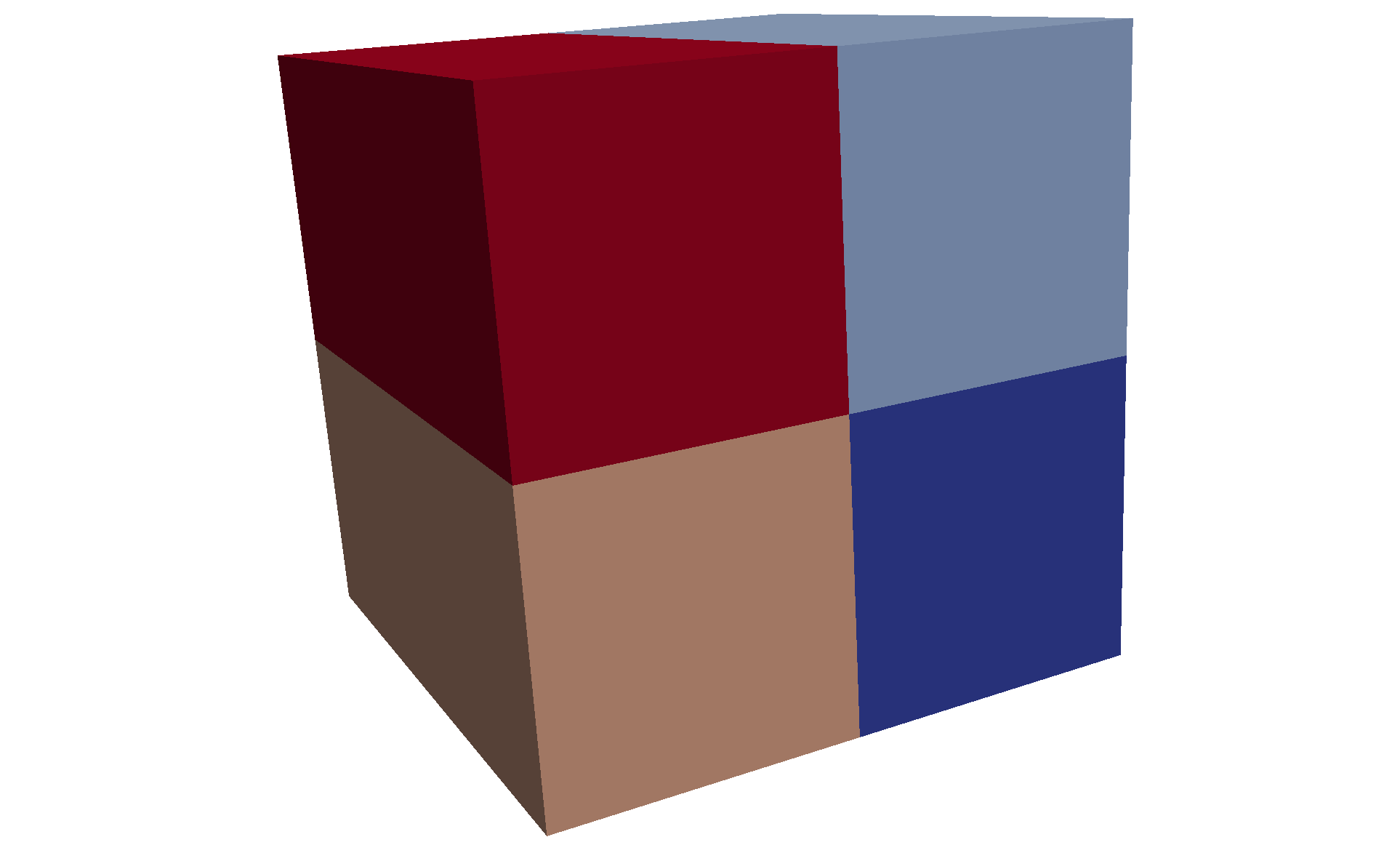}
\includegraphics[width=0.3\textwidth,trim=13cm 0cm 11cm 0cm, clip]{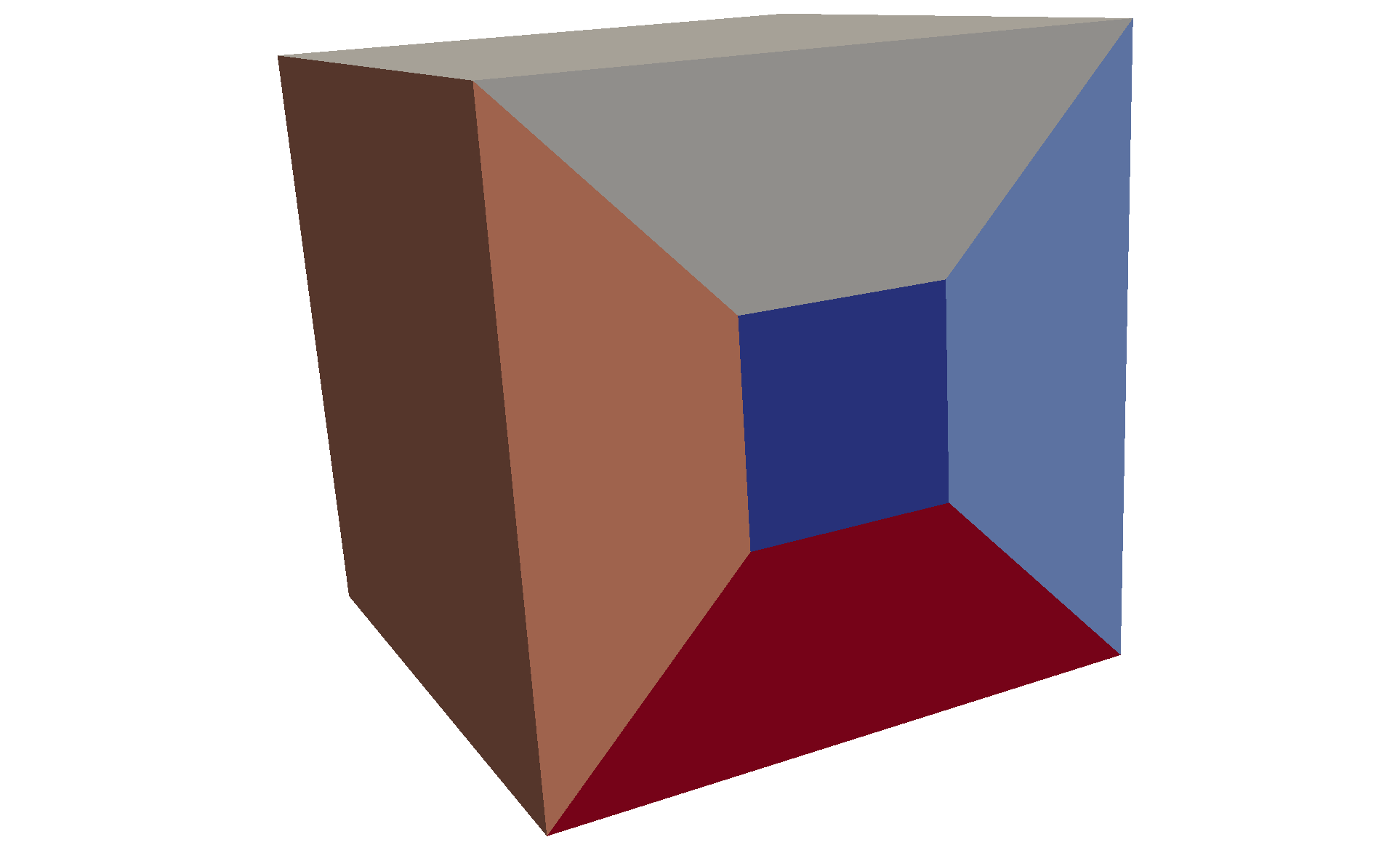}
\caption{Three geometry configurations of the unit cube.}
\label{fig:cube-patches}
\end{figure}

We solve the Maxwell eigenvalue problem on the unit cube, with $\Gamma_{\mathrm{int}}$ one side of the cube, and we consider three different geometric configurations: the first one formed by two patches without internal vertices, the second one with four patches and one internal vertex on $\Gamma_{\mathrm{int}}$, and the third one with five patches and four internal vertices on $\Gamma_{\mathrm{int}}$, see \cref{fig:cube-patches}. The results in \cref{table:zero-eigenfunctions} show the number of zeros obtained in the eigenvalue problem (equal to the dimension of $K_h$), and the expected number of zeros, i.e. the dimension of $\nabla X^0_{h,0}$ (equal to the dimension of $X^0_{h,0}$ due to vanishing boundary conditions), using two different degrees and mesh sizes. We see that in the case of two patches the numbers coincide, while in the others the difference between both is exactly equal to $\# Z_{\Gamma_{\mathrm{int}}}$, independently of the mesh size and the degree of the discrete spaces.

\begin{figure}
\centering
\includegraphics[width=0.3\textwidth,trim=13cm 0cm 11cm 0cm, clip]{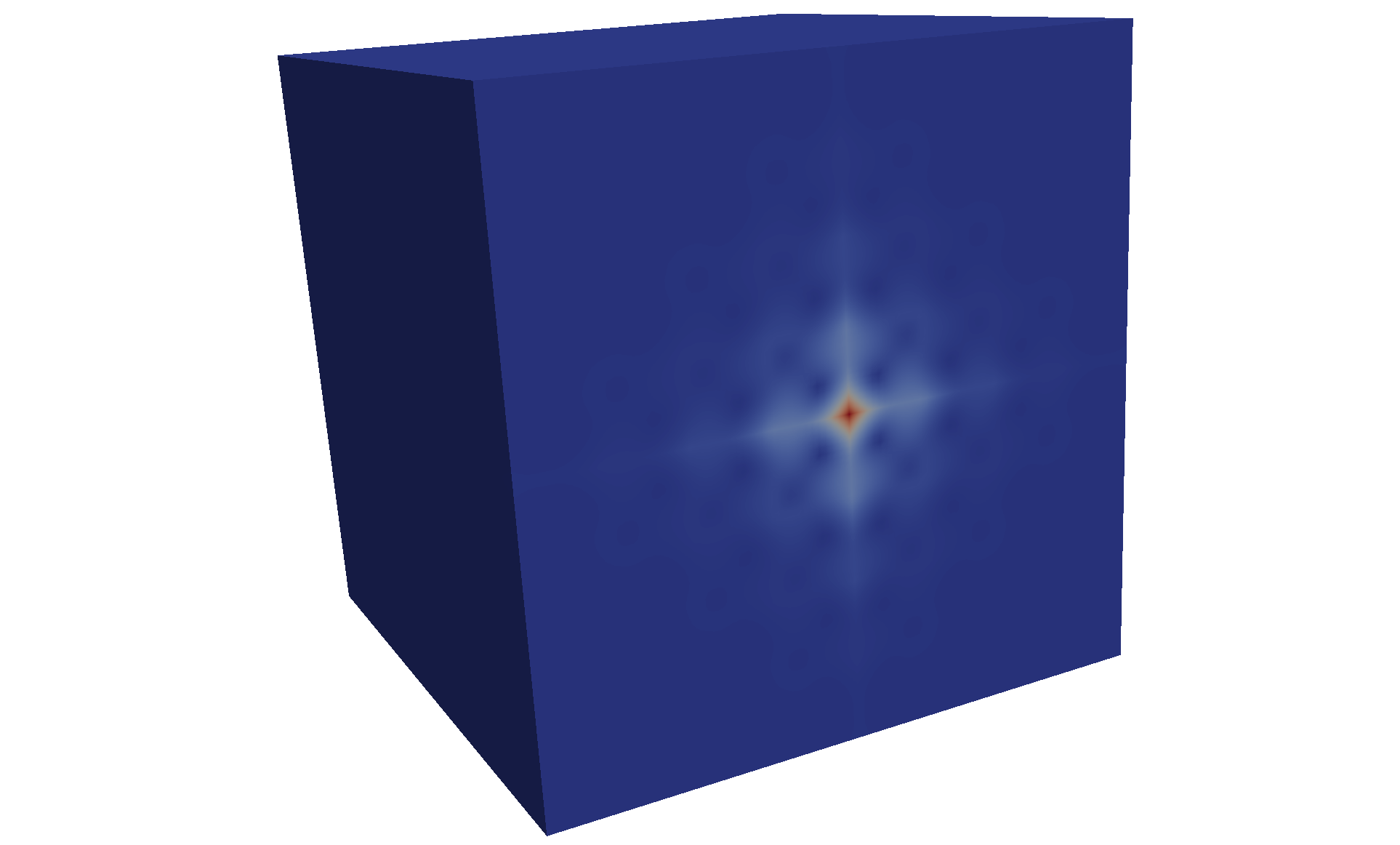}
\includegraphics[width=0.3\textwidth,trim=13cm 0cm 11cm 0cm, clip]{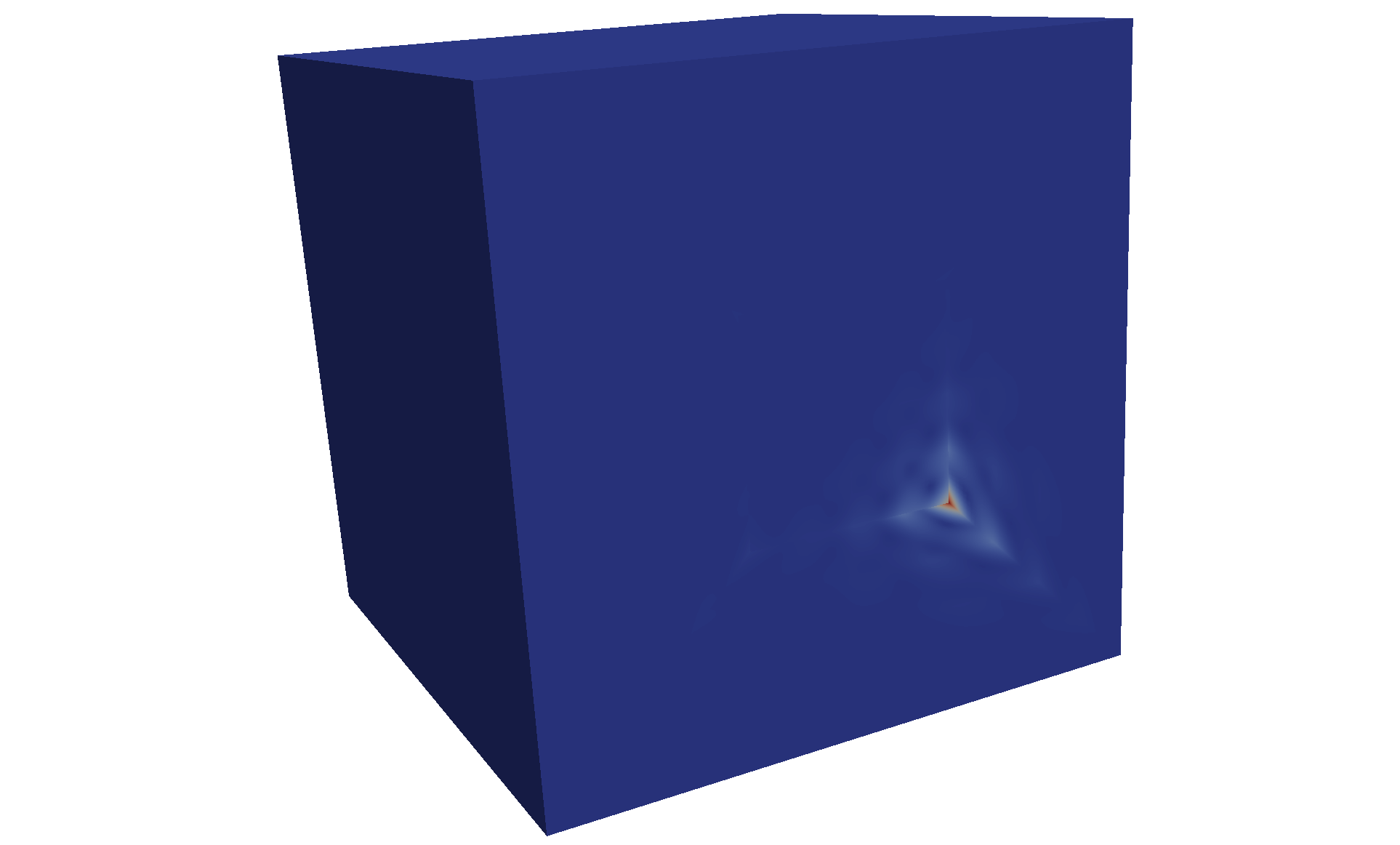}
\includegraphics[width=0.3\textwidth,trim=13cm 0cm 11cm 0cm, clip]{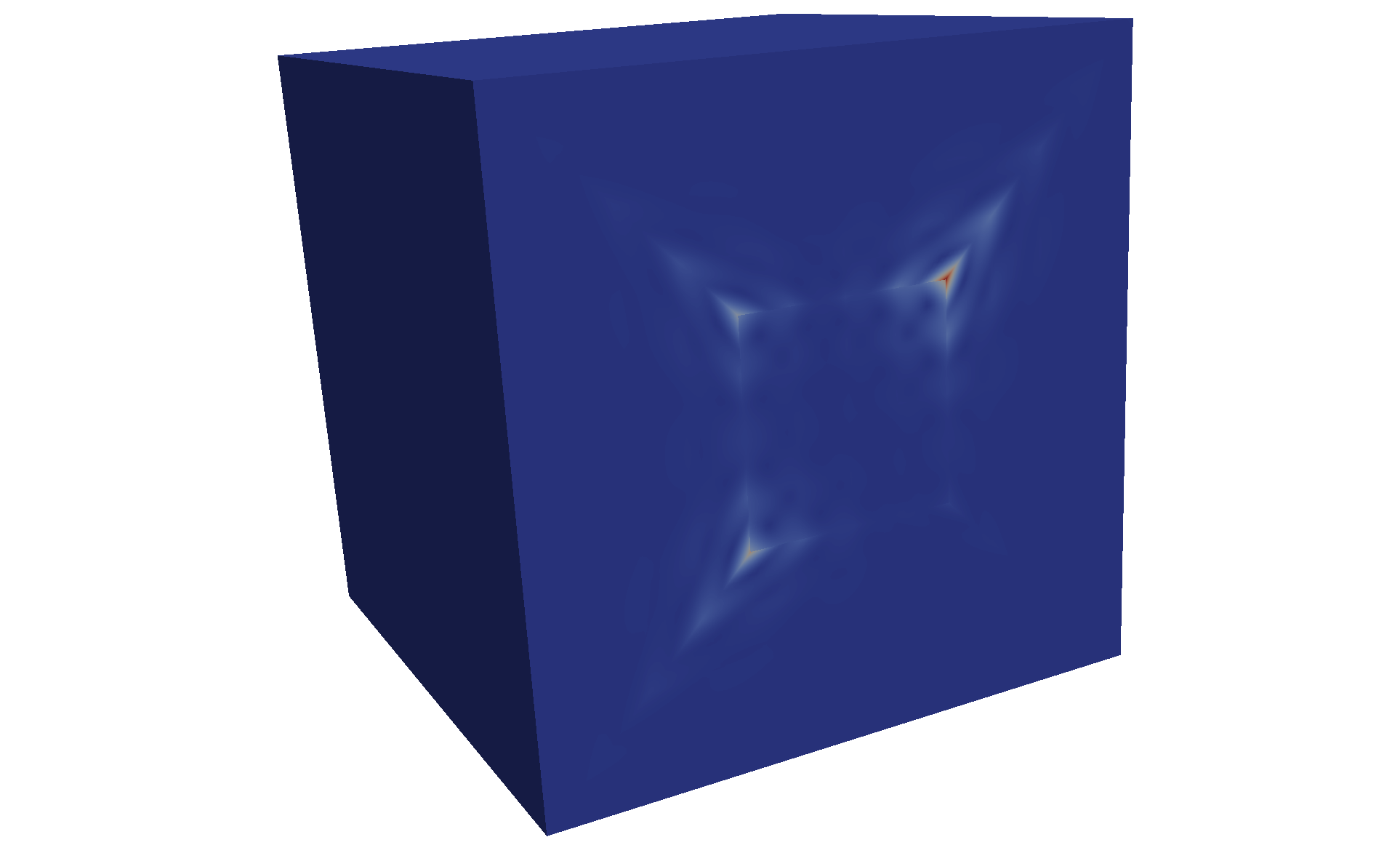}
\caption{Magnitude of some functions in the discrete kernel $K_h$ and not in $\nabla V^0_{h,0}$, for four patches (left) and five patches (center and right).}
\label{fig:cube-functions}
\end{figure}

\begin{table}
\centering
\begin{tabular}{|c|c||c|c|c||c|c|c||c|c|c|}
\hline
\multirow{3}{*}{$p$} & \multirow{3}{*}{$h$} & \multicolumn{3}{|c||}{2 patches} & \multicolumn{3}{|c||}{4 patches} & \multicolumn{3}{|c|}{5 patches} \\
\cline{3-11}
& & \multirow{2}{*}{$\dim X^0_{h,0}$} & \multicolumn{2}{|c||}{$\dim K_h$} & \multirow{2}{*}{$\dim X^0_{h,0}$} & \multicolumn{2}{|c||}{$\dim K_h$} & \multirow{2}{*}{$\dim X^0_{h,0}$} & \multicolumn{2}{|c|}{$\dim K_h$} \\
\cline{4-5} \cline {7-8} \cline{10-11}
& & & $M_h$ & $\widetilde{M}_h$ & & $M_h$ & $\widetilde{M}_h$ & & $M_h$ & $\widetilde{M}_h$ \\
\hline
\multirow{2}{*}{2} 
& $1/2$ & 20 & 20 & 20 & 50 & 51 & 50 & 80 & 84 & 80 \\ \cline{2-11}
& $1/3$ & 63 & 63 & 63 & 147 & 148 & 147 & 219 & 223 & 219\\ \cline{2-11}
\hline
\multirow{2}{*}{3} 
& $1/2$ & 144 & 144 & 144 & 324 & 325 & 324 & 464 & 468 & 464 \\ \cline{2-11}
& $1/3$ & 468 & 468 & 468 & 1014 & 1015 & 1014 & 1392 & 1396 & 1392 \\ \cline{2-11}
\hline
\end{tabular}
\caption{Dimension of $X_{h,0}^0$ and dimension of the discrete kernel $K_h$ using the spaces $M_h$ and $\widetilde{M}_h$, for the three geometry configurations of \cref{fig:cube-patches}.}
\label{table:zero-eigenfunctions}
\end{table}

Using Kikuchi's mixed formulation \cite{Kikuchi} it is easy to remove the zeros associated to $\nabla X^0_{h,0}$. A plot of the eigenfunctions associated to the remaining ones, see \cref{fig:cube-functions}, shows that they oscillate on $\Gamma_{\mathrm{int}}$, with their maxima concentrated on the internal vertices. Since each vertex on $Z_{\Gamma_{\mathrm{int}}}$ exactly coincides with a control point, and has a function associated to it, we simply enrich the space with the gradient of the functions associated to these control points. More precisely, let us denote with $\mathcal{B}_{Z_{\Gamma_{\mathrm{int}}}} \subset X^0_{h}$ the spline basis functions associated to those vertices. We then define the new space for the Lagrange multiplier as 
\[
\widetilde M_h = M_h \oplus G_h,
\]
with $G_h = \mathrm{span} (\nabla_\Gamma \gamma^0_{\Gamma_{\mathrm{int}}} (\mathcal{B}_{Z_{\Gamma_{\mathrm{int}}}}))$, and $M_h$ the space defined above. Notice that each of these functions has support on all the patches adjacent to the corresponding vertex.

\section{Tree-Cotree Gauging}\label{sec:treecotree}
When solving (either directly or iteratively) the linear system resulting from the discretisation of the weak form \eqref{eq:magnetostatic_weakform}, which can be written as
\begin{equation}
\mathbf{K} \mathbf{a} = \mathbf{j},\label{eq:ungauged_system}
\end{equation}
the symmetric system matrix $\mathbf{K}$ will be singular, since it is representing the discrete curl-curl operator. The rank deficiency is due to the lack of uniqueness of the magnetic vector potential $\vecs{A}$, as already mentioned in Section~\ref{sec:magnetostatic}. 
A consistent way to filter out the null-space of the curl operator for lowest-order edge elements is tree-cotree gauging, originally proposed by Albanese and Rubinacci in \cite{Albanese_1988aa}, see also \cite{Bossavit_1998aa,Munteanu_2002aa}. We briefly present the algorithm in the low-order finite element case, and then show how the idea applies seamlessly to the general IGA setting in the single domain case, using the control grid. Finally, we explain the needed modifications for the case of isogeometric mortaring.

\subsection{Basics of Tree-Cotree Gauging in the Finite Element Case}
\subsubsection{Simply Connected Domains}
\begin{figure}[t!]
    \centering
\includegraphics[width=0.4\textwidth,trim=2cm 1cm 2cm 0cm, clip]{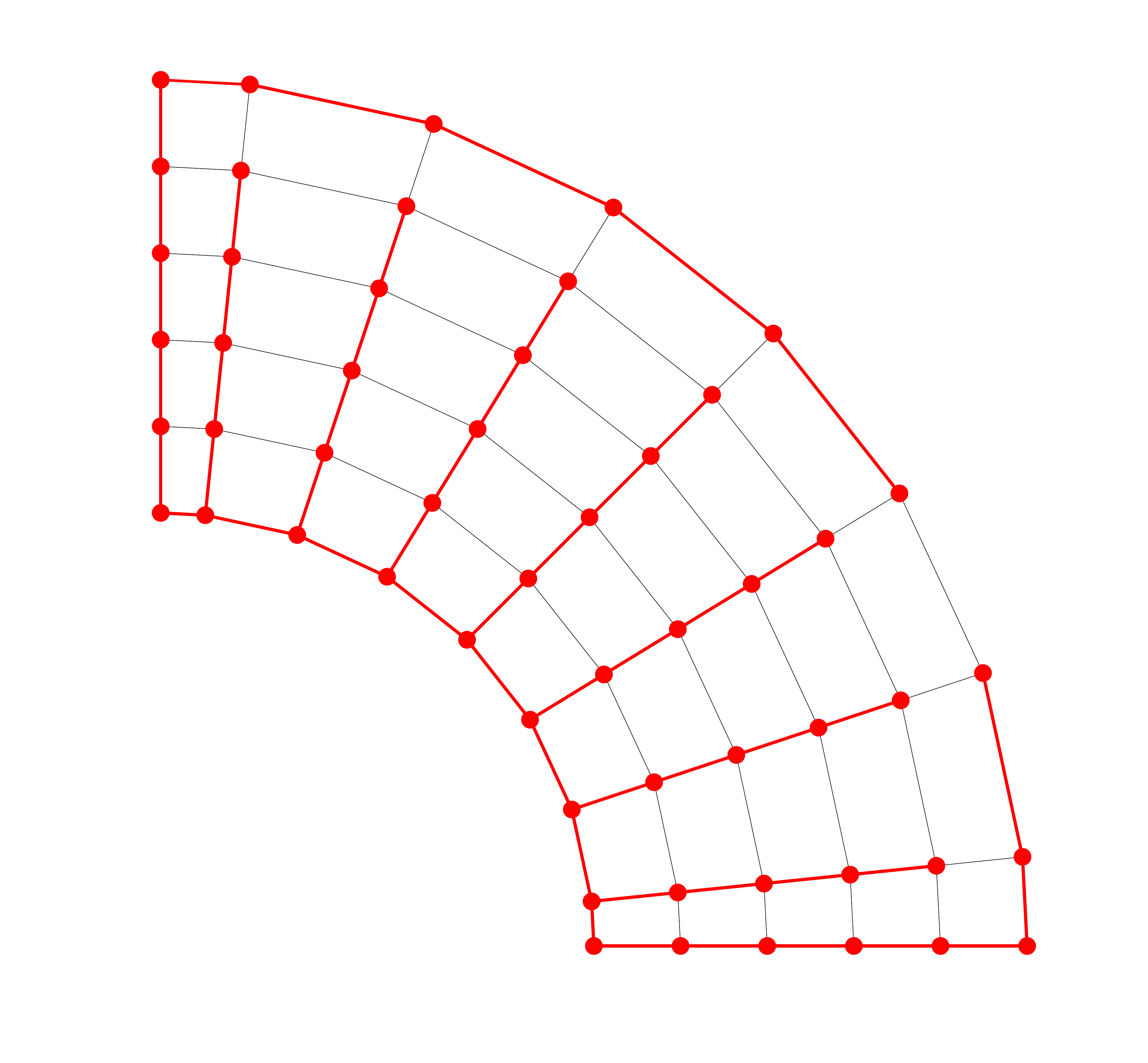}
\caption{The spanning tree built on the graph generated by the vertices and edges of a mesh for the domain in \cref{fig:ctrl_mesh_example}.} \label{fig:tree_example}
\end{figure}

Let us first assume that both the domain $\Omega$ and the Dirichlet boundary $\Gamma_D$ are connected and simply connected. By considering the physical interpretation of degrees of freedom (DoFs) as line integrals of vector fields, and hence their one-to-one correspondence to edges in the mesh, it was shown in \cite{Albanese_1988aa} that setting the value of the DoFs attached to a spanning tree $T$ in the mesh graph, as the one shown in \cref{fig:tree_example} {to zero} yields a reduced system matrix for the complementary set of edges, the cotree $C = E \Setminus T$. This matrix is positive definite. This means seeking a decomposition of \cref{eq:ungauged_system} of the form
\begin{equation}\label{eq:decomposed_system}
    \begin{bmatrix}
        \mathbf{K}_{CC} & \mathbf{K}_{CT}\\
        \mathbf{K}_{TC} & \mathbf{K}_{TT}
        \end{bmatrix}
        \begin{bmatrix}
            \mathbf{a}_{C}\\
            \mathbf{a}_{T}
        \end{bmatrix} =
        \begin{bmatrix}
            \mathbf{j}_{C}\\
            \mathbf{j}_{T}
        \end{bmatrix},
\end{equation}
where the $C$ and $T$ subscripts respectively refer to the cotree and tree edges, such that the reduced system $\mathbf{K}_{CC} \mathbf{a}_C = \mathbf{j}_C$ has a unique solution. Since the solution obtained for the vector potential is of no practical relevance, as only its curl is physically significant, one can set $\mathbf{a}_T = 0$ and solve the system for the basis functions associated to cotree edges. Other choices are discussed for example in \cite{Munteanu_2002aa}.

When building the spanning tree, it is important to correctly deal with Dirichlet boundary conditions, as the corresponding DoFs are imposed strongly and therefore not part of the system \eqref{eq:decomposed_system}. A simple way to do so is by starting the tree from the Dirichlet boundary $\Gamma_D$, and then growing the tree on the rest of the boundary and in the interior, in the following way (see also \cite{dular_discrete_1995}):
\begin{algor} \label{alg:spanning-tree}
    \textbf{Tree-Cotree Decomposition} \begin{enumerate}
\item Build a spanning tree on the subgraph associated to $\Gamma_D$.
\item From the vertices on $\Gamma_D$, grow the tree on $\Gamma_N$.
\item Grow the tree on the remaining graph in the interior of $\Omega$.
\end{enumerate}
\end{algor}
We remark that, since we assumed that both $\Omega$ and $\Gamma_D$ are connected and simply connected, the spanning trees on $\Gamma_D$ and subsequently in $\Omega$ can be grown by simple breadth--first search (BFS) steps \cite{cormen_introduction_2009}. Note also that edges on $\Gamma_D$ which are not part of the grown tree are also not to be added to the cotree, as their associated DoFs are imposed before the solution of the linear system.

\subsubsection{Multiply Connected Domains and Periodic Boundary Conditions}

\begin{figure}[t!]
\begin{subfigure}{0.5\textwidth}
        \includegraphics[width=\textwidth]{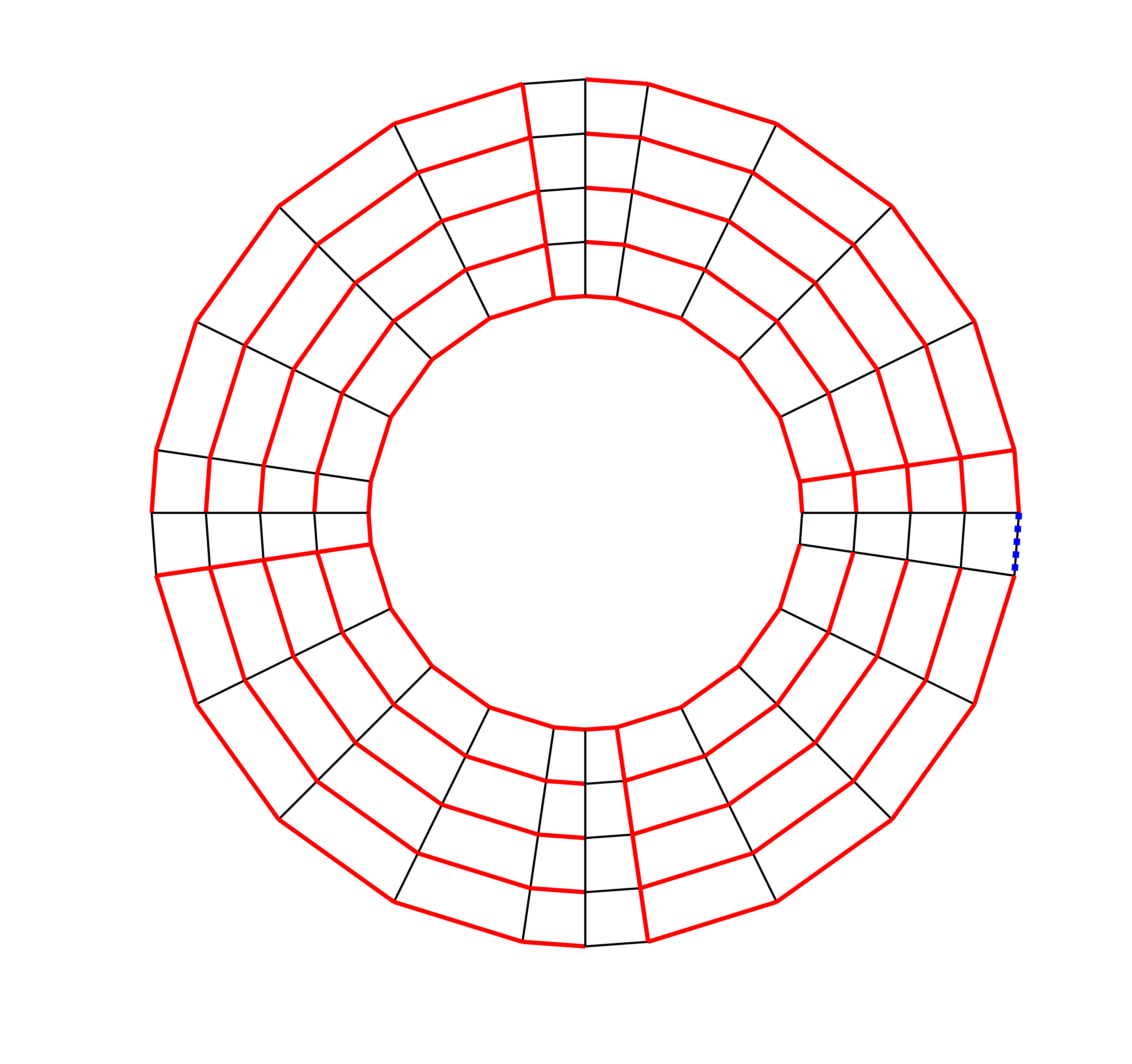}
\caption{Enrichment of the tree for non-trivial topology.}
\label{fig:hole_example}
\end{subfigure}
\begin{subfigure}{0.5\textwidth}
        \includegraphics[width=\textwidth]{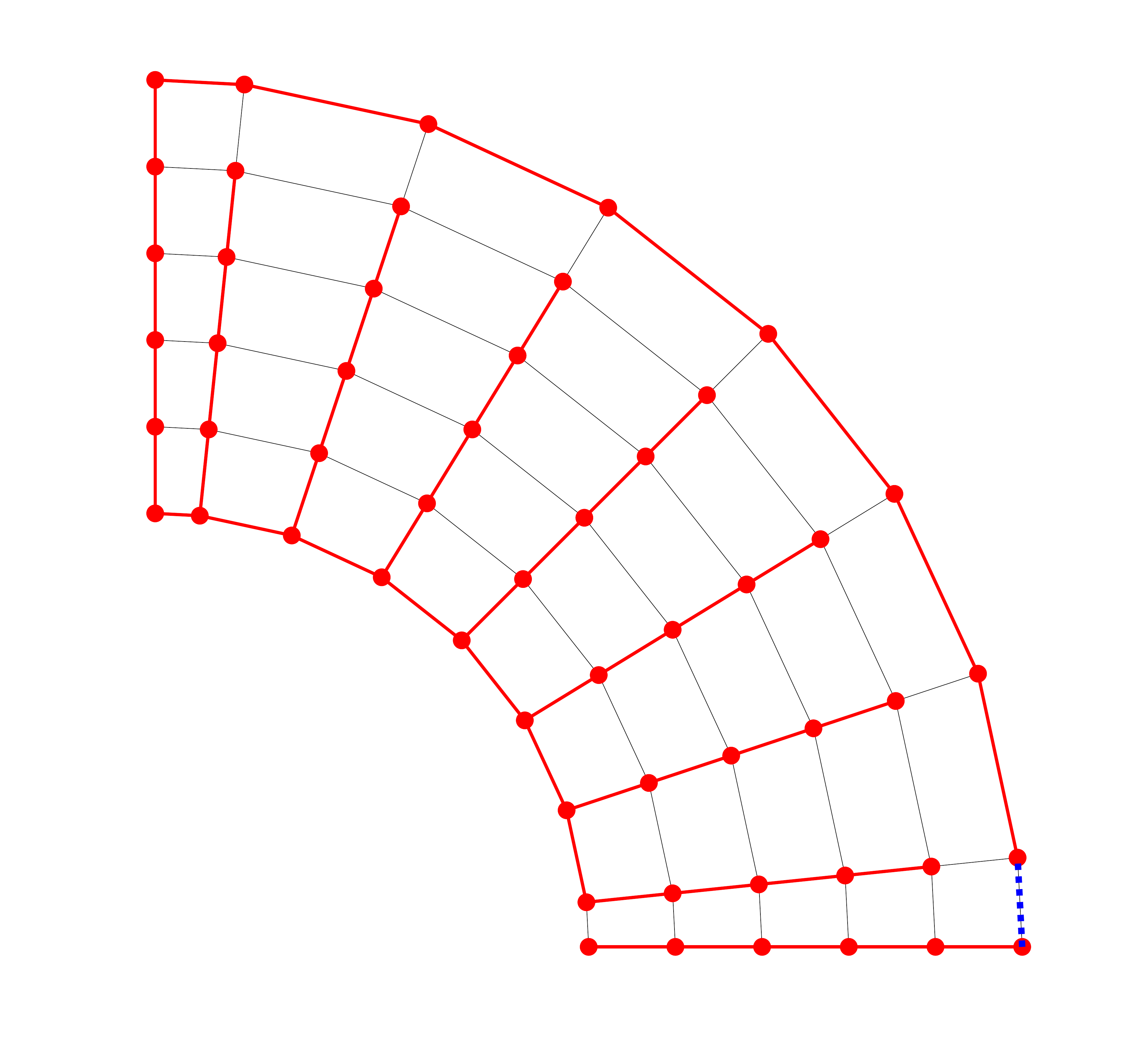}
\caption{Enrichment of the tree for periodic conditions.}
\label{fig:periodic_example}
\end{subfigure}
    \caption{Construction of the tree for non-trivial cases. The spanning tree (in red) is not sufficient for gauging, and it is necessary to add one more function (the one in dotted blue, for example) to the tree.} \label{fig:cohomology_example}
\end{figure}

The spanning tree is sufficient to filter out the discrete kernel when this only consists of gradient functions. This is the case when $\Omega$ and $\Gamma_D$ have trivial topology. In the non-trivial case the discrete kernel of the curl operator will contain, apart from gradients, the space of harmonic fields \cite[Section~1.4]{rodriguez_eddy_2010}. The dimension of this space is equal to the rank of the first cohomology group \cite{munkres_elements_1984}, that is, the number of closed curves that are not a boundary of a surface (up to homology classes). This is exemplified in \cref{fig:hole_example}: the first cohomology group of the annulus is non-trivial, and the spanning tree (in red) is not sufficient to remove all functions with vanishing curl. Adding one more edge to the tree (in dotted blue), removes this additional function, and filters out the discrete kernel. Similar considerations apply when dealing with periodic boundary conditions: in this case the two boundaries where we impose periodicity are identified as a single one, including their nodes and edges, which makes the domain topologically equivalent to an annulus. The spanning tree in \cref{fig:periodic_example}, which takes into account the identification of the periodic boundaries, is not sufficient to gauge the problem, and one more edge (in dotted blue) needs to be added to the tree. In fact, considering the lowest order interpretation as line integrals over edges of the mesh, without that function one could quickly construct a vector potential $\Afield$ such that $\Curl\Afield=0$ everywhere on $\Omega$ yet $\Afield$ is not a gradient.

\begin{remark} \label{rem:cohomology}
In fact, the dimension of the space of harmonic fields is equal to the rank of the first group in relative cohomology with respect to $\Gamma_D$, see again \cite{munkres_elements_1984}, or \cite[Section~5.3]{Bossavit_1998aa} where it is explained as \emph{modulo} $\Gamma_D$. The use of relative cohomology is the reason why the Dirichlet boundary needs to be treated in a different way also in the simply connected case. Relative cohomology must also be used in the case of periodic boundary conditions.
\end{remark}

To overcome the shortcomings of the algorithm based on pure graph theory and not on cohomology, it is necessary to use an extension of the tree--cotree decompositions. The specific algorithm we use is based on existing developments of low order methods and is valid for conventional FEM, the Finite Integration Technique (FIT), the Cell Method (CM), or Discrete Exterior Calculus (DEC). An in-depth description of the procedure goes beyond the scope of the present article and the reader is referred to \cite{dlotko_novel_2013,kapidani_computation_2016,dlotko_topoprocessor_2017,dlotko_lean_2018} for details on their complexity and generality. It suffices to remark here that this additional step needs no \emph{a priori} knowledge of the geometry and topology of $\Omega$, and that the algorithm has linear complexity with respect to the number of DoFs for $\Afield$.

\subsection{Tree-Cotree Gauging in IGA: the Control Grid as a Graph}
In the IGA setting, as soon as $p>1$, the DoFs will lose their physical interpretation as line integrals on edges. Nevertheless, as noted in \cref{sec:iga}, the basis functions of the spline spaces $S_p^k(\Omega)$ in the discrete De Rham sequence \cref{eq:spline_derham} can be identified with $k-$dimensional entities of the control mesh, i.e. with vertices, edges, facets, and cells of the control mesh for $k=0,1,2$ and $3$, respectively.

Thanks to this identification, the spline spaces $S_p^k(\Omega)$ are isomorphic to the lowest order spaces defined in the control mesh. Moreover, as we also mentioned in \cref{sec:iga}, these isomorphisms commute with the derivative. In particular, this implies that a spline function is in the discrete kernel of the curl operator, if and only if its image in the low order edge finite elements in the control mesh is in the kernel of the curl operator for that mesh. Put in simpler terms: filtering out the discrete kernel for splines is equivalent to do it for the lowest order space on the control mesh, and \emph{vice versa}.

It suffices then to apply the algorithm for the tree and cotree construction to the control mesh, and identify the corresponding B-spline functions. Since the algorithm is exactly the same as for finite elements, existing implementations do not need to be changed. Very conveniently, this holds also for different boundary conditions and the topology issues mentioned in the previous section.

It is worth noting that while the application of tree-cotree gauging is automatic for B-splines of any degree, the same is not generally true for high order FEM. In fact, the use of tree-cotree gauging in this setting requires the construction of suitable basis functions as in \cite{Rodriguez_2020aa} or \cite{Schoeberl_2005aa}.

\begin{remark}
The control mesh only plays an auxiliary role to build the spanning tree, and the low order finite element spaces on the control mesh are not needed for any real computation. In fact, the control mesh can be replaced by other meshes with the same topology, for instance, the image of the Greville points through the parametrisation $\mathbf{F}$, also called the Greville mesh in \cite{BSV13}.
\end{remark}

\subsection{Tree-Cotree Gauging for Multi-Patch Domains and Mortaring}
The last step is the combination of tree-cotree gauging with the mortaring technique introduced in \cref{sec:mortar-multipatch}. First of all, we note that the use of tree-cotree gauging for multi-patch domains without mortaring is trivial as long as the parametrisations and the meshes are conforming, in the sense of \cref{sec:multipatch-ungauged}. Under those assumptions, the control meshes from different patches coincide on their interfaces, and the commutative isomorphisms that we defined above between the multi-patch splines spaces and the finite element spaces remain valid using the multi-patch control mesh. Therefore, the algorithm can be applied as before.

\begin{figure}
 \centering
       \tdplotsetmaincoords{0}{0}
      \begin{tikzpicture}[scale=2,tdplot_main_coords,scale = 1.5, xshift=1cm, yshift=0.2cm, rotate around y=-12, rotate around x=15]\def\x{1}
        \def\y{0.5}
        \def\z{1}
        \def\nxl{2}
        \def\nyl{1}
        \def\nzl{2}
        \def\nxu{2}
        \def\nyu{1}
        \def\nzu{2}
        \def\distance{0.3}
        \def\xdistance{2.5}
\def\phirt{60};
        \def\phist{60};
        \def\rirt{1};
        \def\rort{1.5};
        \def\rist{\rort};
        \def\rost{2};
        \def\nphirt{10};
        \def\nphist{7};
        \def\nrrt{5};
        \def\nrst{7};

\begin{scope}[yshift=-\distance cm]
            \draw[blue!20,fill=blue!20] (0,0,0) -- ++(-\x,0,0) -- ++(0,-\y,0) -- ++(\x,0,0) -- cycle;
            \draw[blue!20,fill=blue!20] (0,0,0) -- ++(0,0,-\z) -- ++(0,-\y,0) -- ++(0,0,\z) -- cycle;
            \draw[blue!20,fill=blue!20] (0,0,0) -- ++(-\x,0,0) -- ++(0,0,-\z) -- ++(\x,0,0) -- cycle;
            \draw[blue!20,fill=blue!20] (0,0,0) -- ++(-\x,0,0) -- ++(0,0,-\z) -- ++(\x,0,0) -- cycle;
            {\draw[blue!20,fill=teal!20] (0,-\y,0) -- ++(-\x,0,0) -- ++(0,0,-\z) -- ++(\x,0,0) -- cycle;}
            \node at (0,-\y/2,-\z) [blue,anchor=west] {$\Omega_{1}$};
            \node at (0,-\y,-\z/2) [teal,anchor=west,xshift=5pt] {$\Gamma_{\mathrm{int}}$};
        \end{scope}
\begin{scope}[yshift=-\distance cm]
        \foreach \i in {0,1,...,\nxl}{
            \foreach \j in {0,1,...,\nyl}{
                \foreach \k in {0,1,...,\nzl}{
\draw [gray] (0,-\j*\y/\nyl,-\k*\z/\nzl)--(-\x,-\j*\y/\nyl,-\k*\z/\nzl);
\draw [gray] (-\i*\x/\nxl,0,-\k*\z/\nzl)--(-\i*\x/\nxl,-\y,-\k*\z/\nzl);
\draw [gray] (-\i*\x/\nxl,-\j*\y/\nyl,0)--(-\i*\x/\nxl,-\j*\y/\nyl,-\z);
                }  
            }
        }
{
            \draw [dashed,ultra thick,red] (0,-\y,0)--++(-\x,0,0)--++(0,0,-\z)--++(\x,0,0)--++(0,0,\z/2);
            \draw [dashed,ultra thick,red] (-\x,-\y,-\z/2)--++(\x/2,0,0);
        }
{
            \draw [ultra thick,red] (0,-\y,0)--++(0,\y,0)--++(0,0,-\z)--++(-\x/2,0,0)--++(0,0,\z/2);
            \draw [ultra thick,red] (-\x,-\y,-\z)--++(0,\y,0)--++(0,0,\z)--++(\x/2,0,0);
        }
        \end{scope}
\end{tikzpicture}
  \caption{Visualisation of the tree on the dependent domain $\Omega_{1}$. The tree is shown in red. Dashed lines represent the edges that are built on the interface $\Gamma_{\mathrm{int}}$ but are then removed from the tree.} 
 \label{fig:dependent_tree}
\end{figure}
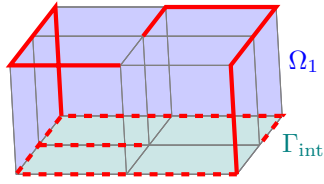

Instead, the validity of tree-cotree gauging when performing mortaring is not obvious, as the \textcolor{black}{dependent} domain $\Omega_1$ and the \textcolor{black}{independent} domain $\Omega_2$ must be treated differently. In fact, one can think that the Lagrange multiplier leaves the solution ``free'' on the interface for the \textcolor{black}{independent} domain, while it ``constrains'' it for the \textcolor{black}{dependent} domain. Therefore, for the \textcolor{black}{independent} domain it is sufficient to use Algorithm~\ref{alg:spanning-tree} replacing $\Gamma_D$ by $\partial \Omega_2 \cap \Gamma_D$, and $\Gamma_N$ by $\partial \Omega_2 \cap (\Gamma_N \cup \Gamma_{\mathrm{int}}) $, i.e. the interface is treated as a Neumann side. Instead, in the \textcolor{black}{dependent} domain the interface will be treated similarly to a Dirichlet side, with the difference that values of DoFs on the interface are not imposed strongly but through the Lagrange multiplier, and therefore they must be added to the cotree as visualised in \cref{fig:dependent_tree}. Therefore, the algorithm for the tree generation in the \textcolor{black}{dependent} subdomain is modified as follows for a simply connected domain:
\begin{algor} \textbf{Modified Tree-Cotree Decomposition} \begin{enumerate}
\item Build a spanning tree on the subgraph associated to $\partial \Omega_1 \cap \Gamma_{\mathrm{int}}$.
\item Grow the tree on the Dirichlet boundary $\partial \Omega_1 \cap (\Gamma_D \Setminus \Gamma_{\mathrm{int}})$.
\item From the vertices on $\partial \Omega_1 \cap (\Gamma_D \cup \Gamma_{\mathrm{int}})$, grow the tree on $\partial \Omega_1 \cap \Gamma_N$.
\item Grow the tree on the remaining graph in the interior of $\Omega_1$.
\item Restore all edges supported on $\Gamma_{\mathrm{int}}$ as cotree DoFs.
\end{enumerate}
\end{algor}

Finally, restricting directly the linear system to the cotree DoFs, the linear system that we have to solve in the mortaring case takes the form
\begin{equation} \label{eq:matrix-gauged}
    \begin{bmatrix}
        \mathbf{K}^1_{CC} & & (\mathbf{G}^1_C)^\top\\
        & \mathbf{K}^2_{CC} & (\mathbf{G}^2_C)^\top\\
        \mathbf{G}^1_C & \mathbf{G}^2_C
        \end{bmatrix}
        \begin{bmatrix}
            \mathbf{a}^1_{C}\\
            \mathbf{a}^2_{C}\\
            \boldsymbol{\mu}
        \end{bmatrix} =
        \begin{bmatrix}
            \mathbf{j}^1_{C}\\
            \mathbf{j}^2_{C}\\
            \mathbf{0}
        \end{bmatrix},
\end{equation}
where the numerical superscripts refer to each subdomain, and the $C$ subindex refers to the \emph{cotree} DoFs in the corresponding subdomain.
\textcolor{black}{Note also that edges on $\partial \Omega_1 \cap \Gamma_D$ which are not part of the grown tree are also not to part of the final cotree, as their associated DoFs are imposed before the solution of the linear system.}
\begin{remark}
We have restricted ourselves to two domains just for simplicity, but mortaring can be applied to a situation with more domains, as reported in \cite{Buffa_2020aa}. Accordingly, tree-cotree gauging can also be used in this more general context, and we note that in this case the tree-cotree decomposition for the different subdomains can be carried out in parallel.
\end{remark}

\section{Numerical Results}\label{sec:numerics}
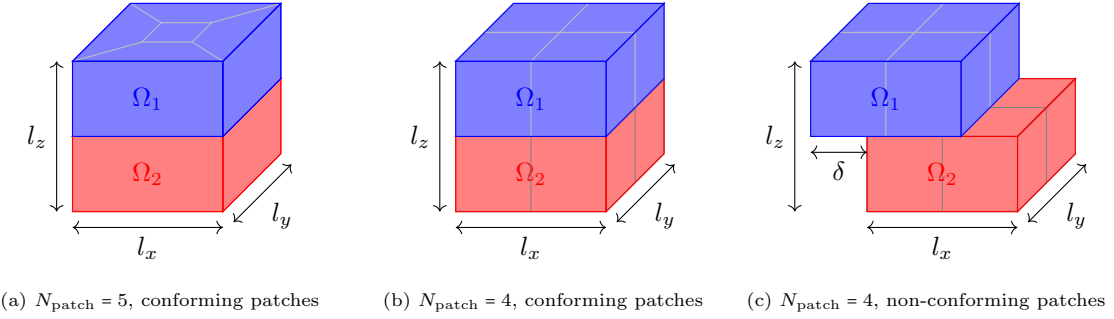
\begin{figure}
\centering
\begin{subfigure}{0.3\textwidth}
    \centering 
\begin{tikzpicture}[scale = 2, xshift=1cm, yshift=0.2cm]\def\x{1}
        \def\y{0.5}
        \def\z{1}
        \def\nxl{2}
        \def\nyl{1}
        \def\nzl{2}
        \def\nxu{2}
        \def\nyu{1}
        \def\nzu{2}
        \def\distance{0}
        \def\xdistance{2.5}
\begin{scope}[yshift=-\distance cm]
            \draw[red,fill=red!50] (0,0,0) -- ++(-\x,0,0) -- ++(0,-\y,0) -- ++(\x,0,0) -- cycle;
			\node[red] at (-\x/2,-\y/2,0) {$\Omega_{2}$};
            \draw[red,fill=red!50] (0,0,0) -- ++(0,0,-\z) -- ++(0,-\y,0) -- ++(0,0,\z) -- cycle;
            \draw[red,fill=red!50] (0,0,0) -- ++(-\x,0,0) -- ++(0,0,-\z) -- ++(\x,0,0) -- cycle;

\draw[red] (0,0,0) -- ++(-\x,0,0) -- ++(0,-\y,0) -- ++(\x,0,0) -- cycle;
            \draw[red] (0,0,0) -- ++(0,0,-\z) -- ++(0,-\y,0) -- ++(0,0,\z) -- cycle;
            \draw[red] (0,0,0) -- ++(-\x,0,0) -- ++(0,0,-\z) -- ++(\x,0,0) -- cycle;

            \draw[<->,yshift=-3pt] (-\x,-\y,0) --++ (\x,0,0) node [pos=0.5,anchor=north] {$l_{x}$}; 
		  \draw[<->,xshift=-3pt] (-\x,-\y,0) --++ (0,2*\y,0) node [pos=0.5,anchor=east] {$l_{z}$}; 
		  \draw[<->,xshift=2pt,yshift=-2pt] (0,-\y,0) --++ (0,0,-\z) node [pos=0.5,anchor=north west] {$l_{y}$}; 
        \end{scope}

\begin{scope}[yshift=\distance cm +\y cm]
            \draw[blue,fill=blue!50] (0,0,0) -- ++(-\x,0,0) -- ++(0,-\y,0) -- ++(\x,0,0) -- cycle;
            \draw[blue,fill=blue!50] (0,0,0) -- ++(0,0,-\z) -- ++(0,-\y,0) -- ++(0,0,\z) -- cycle;
            \draw[blue,fill=blue!50] (0,0,0) -- ++(-\x,0,0) -- ++(0,0,-\z) -- ++(\x,0,0) -- cycle;

            \draw[lightgray] (-\x/3,0,-\z/3) -- ++(-\x/3,0,0) -- ++(0,0,-\z/3) -- ++(\x/3,0,0) -- cycle;
            \draw[lightgray] (-\x/3,0,-\z/3) -- ++(\x/3,0,\z/3);
            \draw[lightgray] (-\x/3,0,-2*\z/3) -- ++(\x/3,0,-\z/3);
            \draw[lightgray] (-2*\x/3,0,-2*\z/3) -- ++(-\x/3,0,-\z/3);
            \draw[lightgray] (-2*\x/3,0,-\z/3) -- ++(-\x/3,0,\z/3);

\draw[blue] (0,0,0) -- ++(-\x,0,0) -- ++(0,-\y,0) -- ++(\x,0,0) -- cycle;
            \draw[blue] (0,0,0) -- ++(0,0,-\z) -- ++(0,-\y,0) -- ++(0,0,\z) -- cycle;
            \draw[blue] (0,0,0) -- ++(-\x,0,0) -- ++(0,0,-\z) -- ++(\x,0,0) -- cycle;

			\node [blue] at (-\x/2,-\y/2,0) {$\Omega_{1}$};
        
\end{scope}
        
    \end{tikzpicture}
  \caption{$N_{\mathrm{patch}}=5$, conforming patches}
\label{fig:problemdomain_cube_fivepatch}
 \end{subfigure}
 \begin{subfigure}{0.3\textwidth}
    \centering 
\begin{tikzpicture}[scale = 2, xshift=1cm, yshift=0.2cm]\def\x{1}
        \def\y{0.5}
        \def\z{1}
        \def\nxl{2}
        \def\nyl{1}
        \def\nzl{2}
        \def\nxu{2}
        \def\nyu{1}
        \def\nzu{2}
        \def\distance{0}
        \def\xdistance{2.5}
\begin{scope}[yshift=-\distance cm]
            \draw[red,fill=red!50] (0,0,0) -- ++(-\x,0,0) -- ++(0,-\y,0) -- ++(\x,0,0) -- cycle;
			\node[red] at (-\x/2,-\y/2,0) {$\Omega_{2}$};
            \draw[red,fill=red!50] (0,0,0) -- ++(0,0,-\z) -- ++(0,-\y,0) -- ++(0,0,\z) -- cycle;
            \draw[red,fill=red!50] (0,0,0) -- ++(-\x,0,0) -- ++(0,0,-\z) -- ++(\x,0,0) -- cycle;

\foreach \i in {0,1,...,\nyl}{
                \draw [gray] (0,-\i*\y/\nyl,0)--(-\x,-\i*\y/\nyl,0);
            }
\foreach \i in {0,1,...,\nzl}{
                \draw [gray] (0,0,-\i*\z/\nzl)--(-\x,0,-\i*\z/\nzl);
            }

\foreach \i in {0,1,...,\nxl}{
                \draw [gray] (-\i*\x/\nxl,0,0)--(-\i*\x/\nxl,-\y,0);
            }
\foreach \i in {0,1,...,\nzl}{
                \draw [gray] (0,0,-\i*\z/\nzl)--(0,-\y,-\i*\z/\nzl);
            }

\foreach \i in {0,1,...,\nxl}{
                \draw [gray] (-\i*\x/\nxl,0,0)--(-\i*\x/\nxl,0,-\z);
            }
\foreach \i in {0,1,...,\nyl}{
                \draw [gray] (0,-\i*\y/\nyl,0)--(0,-\i*\y/\nyl,-\z);
            }
            
\draw[red] (0,0,0) -- ++(-\x,0,0) -- ++(0,-\y,0) -- ++(\x,0,0) -- cycle;
            \draw[red] (0,0,0) -- ++(0,0,-\z) -- ++(0,-\y,0) -- ++(0,0,\z) -- cycle;
            \draw[red] (0,0,0) -- ++(-\x,0,0) -- ++(0,0,-\z) -- ++(\x,0,0) -- cycle;

            \draw[<->,yshift=-3pt] (-\x,-\y,0) --++ (\x,0,0) node [pos=0.5,anchor=north] {$l_{x}$}; 
		  \draw[<->,xshift=-3pt] (-\x,-\y,0) --++ (0,2*\y,0) node [pos=0.5,anchor=east] {$l_{z}$}; 
		  \draw[<->,xshift=2pt,yshift=-2pt] (0,-\y,0) --++ (0,0,-\z) node [pos=0.5,anchor=north west] {$l_{y}$}; 
        \end{scope}

\begin{scope}[yshift=\distance cm +\y cm]
            \draw[blue,fill=blue!50] (0,0,0) -- ++(-\x,0,0) -- ++(0,-\y,0) -- ++(\x,0,0) -- cycle;
            \draw[blue,fill=blue!50] (0,0,0) -- ++(0,0,-\z) -- ++(0,-\y,0) -- ++(0,0,\z) -- cycle;
            \draw[blue,fill=blue!50] (0,0,0) -- ++(-\x,0,0) -- ++(0,0,-\z) -- ++(\x,0,0) -- cycle;

\foreach \i in {0,1,...,\nyu}{
                \draw [lightgray] (0,-\i*\y/\nyu,0)--(-\x,-\i*\y/\nyu,0);
            }
\foreach \i in {0,1,...,\nzu}{
                \draw [lightgray] (0,0,-\i*\z/\nzu)--(-\x,0,-\i*\z/\nzu);
            }

\foreach \i in {0,1,...,\nxu}{
                \draw [lightgray] (-\i*\x/\nxu,0,0)--(-\i*\x/\nxu,-\y,0);
            }
\foreach \i in {0,1,...,\nzu}{
                \draw [lightgray] (0,0,-\i*\z/\nzu)--(0,-\y,-\i*\z/\nzu);
            }

\foreach \i in {0,1,...,\nxu}{
                \draw [lightgray] (-\i*\x/\nxu,0,0)--(-\i*\x/\nxu,0,-\z);
            }
\foreach \i in {0,1,...,\nyu}{
                \draw [lightgray] (0,-\i*\y/\nyu,0)--(0,-\i*\y/\nyu,-\z);
            }

\draw[blue] (0,0,0) -- ++(-\x,0,0) -- ++(0,-\y,0) -- ++(\x,0,0) -- cycle;
            \draw[blue] (0,0,0) -- ++(0,0,-\z) -- ++(0,-\y,0) -- ++(0,0,\z) -- cycle;
            \draw[blue] (0,0,0) -- ++(-\x,0,0) -- ++(0,0,-\z) -- ++(\x,0,0) -- cycle;
           
			\node [blue] at (-\x/2,-\y/2,0) {$\Omega_{1}$};
        
\end{scope}
        
    \end{tikzpicture}
  \caption{$N_{\mathrm{patch}}=4$, conforming patches}
\label{fig:problemdomain_cube_conf}
 \end{subfigure}
\begin{subfigure}{0.3\textwidth}
    \centering
\begin{tikzpicture}[scale = 2, xshift=1cm, yshift=0.2cm]\def\x{1}
        \def\y{0.5}
        \def\z{1}
        \def\nxl{2}
        \def\nyl{1}
        \def\nzl{2}
        \def\nxu{2}
        \def\nyu{1}
        \def\nzu{2}
        \def\distance{0}
        \def\xdistance{2.5}
\begin{scope}[yshift=-\distance cm, xshift=   0.375\x cm]
            \draw[red,fill=red!50] (0,0,0) -- ++(-\x,0,0) -- ++(0,-\y,0) -- ++(\x,0,0) -- cycle;
            \draw[red,fill=red!50] (0,0,0) -- ++(0,0,-\z) -- ++(0,-\y,0) -- ++(0,0,\z) -- cycle;
            \draw[red,fill=red!50] (0,0,0) -- ++(-\x,0,0) -- ++(0,0,-\z) -- ++(\x,0,0) -- cycle;

\foreach \i in {0,1,...,\nyl}{
                \draw [gray] (0,-\i*\y/\nyl,0)--(-\x,-\i*\y/\nyl,0);
            }
\foreach \i in {0,1,...,\nzl}{
                \draw [gray] (0,0,-\i*\z/\nzl)--(-\x,0,-\i*\z/\nzl);
            }

\foreach \i in {0,1,...,\nxl}{
                \draw [gray] (-\i*\x/\nxl,0,0)--(-\i*\x/\nxl,-\y,0);
            }
\foreach \i in {0,1,...,\nzl}{
                \draw [gray] (0,0,-\i*\z/\nzl)--(0,-\y,-\i*\z/\nzl);
            }

\foreach \i in {0,1,...,\nxl}{
                \draw [gray] (-\i*\x/\nxl,0,0)--(-\i*\x/\nxl,0,-\z);
            }
\foreach \i in {0,1,...,\nyl}{
                \draw [gray] (0,-\i*\y/\nyl,0)--(0,-\i*\y/\nyl,-\z);
            }

\draw[red] (0,0,0) -- ++(-\x,0,0) -- ++(0,-\y,0) -- ++(\x,0,0) -- cycle;
            \draw[red] (0,0,0) -- ++(0,0,-\z) -- ++(0,-\y,0) -- ++(0,0,\z) -- cycle;
            \draw[red] (0,0,0) -- ++(-\x,0,0) -- ++(0,0,-\z) -- ++(\x,0,0) -- cycle;

		  \draw[<->,yshift=-3pt] (-\x,-\y,0) --++ (\x,0,0) node [pos=0.5,anchor=north] {$l_{x}$}; 
		  \draw[<->,xshift=-3pt- 0.375 cm] (-\x,-\y,0) --++ (0,2*\y,0) node [pos=0.5,anchor=east] {$l_{z}$}; 
		  \draw[<->,xshift=2pt,yshift=-2pt] (0,-\y,0) --++ (0,0,-\z) node [pos=0.5,anchor=north west] {$l_{y}$}; 
			\node[red] at (-\x/2,-\y/2,0) {$\Omega_{2}$};

        \end{scope}

\begin{scope}[yshift=\distance cm +\y cm]
            \draw[blue,fill=blue!50] (0,0,0) -- ++(-\x,0,0) -- ++(0,-\y,0) -- ++(\x,0,0) -- cycle;
            \draw[blue,fill=blue!50] (0,0,0) -- ++(0,0,-\z) -- ++(0,-\y,0) -- ++(0,0,\z) -- cycle;
            \draw[blue,fill=blue!50] (0,0,0) -- ++(-\x,0,0) -- ++(0,0,-\z) -- ++(\x,0,0) -- cycle;
			  \draw[<->,yshift=-3pt] (-\x,-\y,0) --++ (0.375\x,0,0) node [pos=0.5,anchor=north] {$\delta$}; 
		 
\foreach \i in {0,1,...,\nyu}{
                \draw [lightgray] (0,-\i*\y/\nyu,0)--(-\x,-\i*\y/\nyu,0);
            }
\foreach \i in {0,1,...,\nzu}{
                \draw [lightgray] (0,0,-\i*\z/\nzu)--(-\x,0,-\i*\z/\nzu);
            }

\foreach \i in {0,1,...,\nxu}{
                \draw [lightgray] (-\i*\x/\nxu,0,0)--(-\i*\x/\nxu,-\y,0);
            }
\foreach \i in {0,1,...,\nzu}{
                \draw [lightgray] (0,0,-\i*\z/\nzu)--(0,-\y,-\i*\z/\nzu);
            }

\foreach \i in {0,1,...,\nxu}{
                \draw [lightgray] (-\i*\x/\nxu,0,0)--(-\i*\x/\nxu,0,-\z);
            }
\foreach \i in {0,1,...,\nyu}{
                \draw [lightgray] (0,-\i*\y/\nyu,0)--(0,-\i*\y/\nyu,-\z);
            }

\draw[blue] (0,0,0) -- ++(-\x,0,0) -- ++(0,-\y,0) -- ++(\x,0,0) -- cycle;
            \draw[blue] (0,0,0) -- ++(0,0,-\z) -- ++(0,-\y,0) -- ++(0,0,\z) -- cycle;
            \draw[blue] (0,0,0) -- ++(-\x,0,0) -- ++(0,0,-\z) -- ++(\x,0,0) -- cycle;

			\node [blue] at (-\x/2,-\y/2,0) {$\Omega_{1}$};
\end{scope}
        
    \end{tikzpicture}
  \caption{$N_{\mathrm{patch}}=4$, non-conforming patches}
\label{fig:problemdomain_cube_nonconf}
 \end{subfigure}
\caption{Domain of the test problem with $N_{\mathrm{patch}}$ patches both in $\Omega_{1}$ and $\Omega_{2}$ with conforming patches (a,b) and with non-conforming patches (c) where $\Omega_{2}$ is shifted by $\delta$ in $x$-direction. The parameters of the domain are $l_{x}=\pi$, $l_{y}=\pi$, $l_{z}=\pi$.}
\label{fig:problemdomain_cube}
\end{figure}

As numerical test cases we solve the Maxwell eigenvalue problem \cref{eq:eigvalprob} and the magnetostatic source problem \cref{eq:mortar1,eq:mortar2,eq:mortarcontinuity} on the domain shown in \cref{fig:problemdomain_cube}. All numerical examples are computed using the open-source software GeoPDEs, \cite{Vazquez_2016aa}.

\subsection{Inf-Sup Condition}
\begin{figure}
 \centering
\begin{tikzpicture}[]
           \begin{loglogaxis}[width=0.6\columnwidth, height=0.4\columnwidth, ylabel={$\beta_{\mathrm{inf-sup}}$},xlabel={$1/h$}, ylabel near ticks,tick label style={font=\footnotesize},  legend style={font=\footnotesize, }, label style={font=\footnotesize},every x tick scale label/.style={at={(1,0)},anchor=north,yshift=-5pt,inner sep=0pt},legend pos=outer north east]
                \addplot [blue, mark=triangle*, mark options={solid}] table [x index=0, y index=1, col sep=comma] {images/beta_M_nsub_slavedeg_2_3_4.csv};
                \addplot [red, mark=*, mark options={solid}] table [x index=0, y index=2, col sep=comma] {images/beta_M_nsub_slavedeg_2_3_4.csv};
                \addplot [green, mark=square*, mark options={solid}] table [x index=0, y index=3, col sep=comma] {images/beta_M_nsub_slavedeg_2_3_4.csv};  
                \addplot [dashed, blue, mark=triangle*, mark options={solid}] table [x index=0, y index=1, col sep=comma] {images/beta_M_tilde_nsub_slavedeg_2_3_4.csv};
                \addplot [dashed, red, mark=*, mark options={solid}] table [x index=0, y index=2, col sep=comma] {images/beta_M_tilde_nsub_slavedeg_2_3_4.csv};
                \addplot [dashed, green, mark=square*, mark options={solid}] table [x index=0, y index=3, col sep=comma] {images/beta_M_tilde_nsub_slavedeg_2_3_4.csv};    
                                \addplot [magenta, mark=triangle*, mark options={solid}] table [x index=0, y index=1, col sep=comma] {images/beta_fivepatch_M_nsub_slavedeg_2_3_4.csv};
                \addplot [orange, mark=*, mark options={solid}] table [x index=0, y index=2, col sep=comma] {images/beta_fivepatch_M_nsub_slavedeg_2_3_4.csv};
                \addplot [teal, mark=square*, mark options={solid}] table [x index=0, y index=3, col sep=comma] {images/beta_fivepatch_M_nsub_slavedeg_2_3_4.csv};  
                \addplot [dashed, magenta, mark=triangle*, mark options={solid}] table [x index=0, y index=1, col sep=comma] {images/beta_fivepatch_M_tilde_nsub_slavedeg_2_3_4.csv};
                \addplot [dashed, orange, mark=*, mark options={solid}] table [x index=0, y index=2, col sep=comma] {images/beta_fivepatch_M_tilde_nsub_slavedeg_2_3_4.csv};
                \addplot [dashed, teal, mark=square*, mark options={solid}] table [x index=0, y index=3, col sep=comma] {images/beta_fivepatch_M_tilde_nsub_slavedeg_2_3_4.csv};    
                \legend{$N_{\mathrm{patch}}=4, p=2, M_h$\\$N_{\mathrm{patch}}=4, p=3, M_h$\\$N_{\mathrm{patch}}=4, p=4, M_h$\\$N_{\mathrm{patch}}=4, p=2, \widetilde{M}_h$\\$N_{\mathrm{patch}}=4, p=3, \widetilde{M}_h$\\$N_{\mathrm{patch}}=4, p=4, \widetilde{M}_h$\\$N_{\mathrm{patch}}=5, p=2, M_h$\\$N_{\mathrm{patch}}=5, p=3, M_h$\\$N_{\mathrm{patch}}=5, p=4, M_h$\\$N_{\mathrm{patch}}=5, p=2, \widetilde{M}_h$\\$N_{\mathrm{patch}}=5, p=3, \widetilde{M}_h$\\$N_{\mathrm{patch}}=5, p=4, \widetilde{M}_h$\\
                }
            \end{loglogaxis}
\end{tikzpicture} 
\caption{Inf-sup constant $\beta_{\mathrm{inf-sup}}$ for the multipatch mortaring of the unit cube shown in \cref{fig:problemdomain_cube} with $N_{\mathrm{patch}}$ patches in each subdomain using $M_h$ and $\widetilde{M}_h$ as space for the Lagrange multipliers for different degree $p$.}
 \label{fig:infsup}
\end{figure}
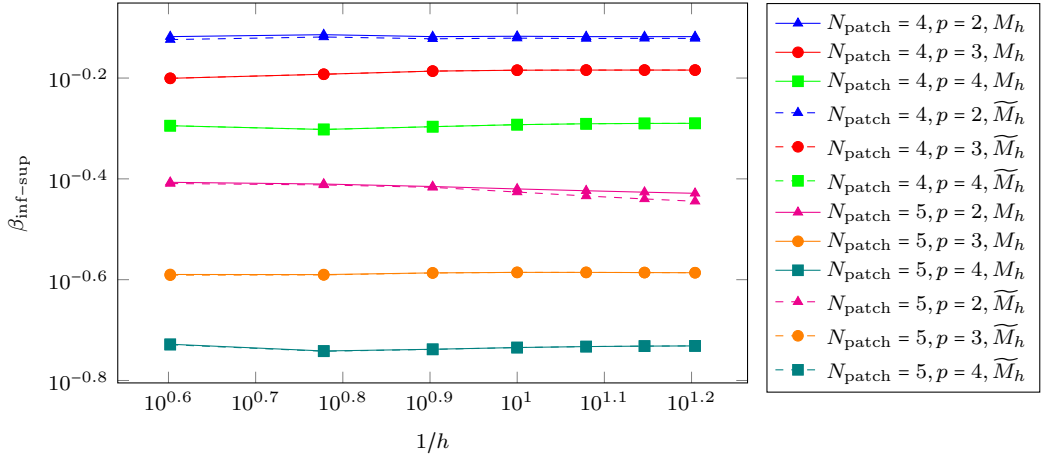

The stability of the saddle-point problem has been investigated in \cite{Buffa_2020aa}, where it has been numerically shown to be stable for $M_h$. To investigate the stability of the problem with the modified Lagrange multiplier space $\widetilde{M}_h$, the inf-sup constant $\beta_{\mathrm{inf-sup}}$ is evaluated numerically \cite{Chapelle_1993aa} for different mesh refinement levels with mesh size $h$ in the parametric domain. The eigenvalue problem \cref{eq:eigvalprob} is solved on the domain shown in \cref{fig:problemdomain_cube_conf} and \cref{fig:problemdomain_cube_fivepatch} with Dirichlet boundary conditions on all boundaries. The results are shown in \cref{fig:infsup}. The inf-sup constant is shown to be bounded away from zero confirming the stability of the problem.
\subsection{Maxwell Eigenvalue Problem}
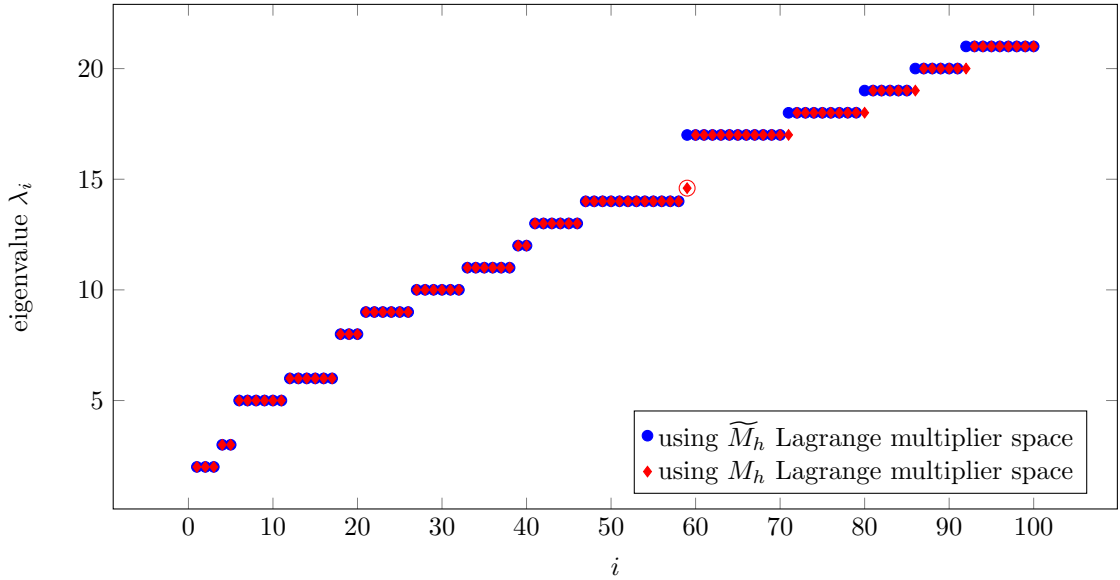
\begin{figure}
\centering
\begin{tikzpicture}
	\begin{axis}[width=0.9\linewidth,height=0.5\linewidth,
		xlabel=$i$,
		ylabel=eigenvalue $\lambda_{i}$,
		legend pos=south east]
        \addplot [blue, mark=*, mark options={solid}, only marks] table [x index=0, y index=3, col sep=comma] {images/eigvals_correct_spurious.csv};
        \addplot [red, mark=diamond*, mark options={solid}, only marks] table [x index=0, y index=4, col sep=comma] {images/eigvals_correct_spurious.csv};
        \addplot [red, mark=o, mark options={solid}, only marks,mark size=3pt] coordinates {(59, 4*3.6495)};
        \legend{using $\widetilde{M}_h$ Lagrange multiplier space , using $M_h$ Lagrange multiplier space }
	\end{axis}
\end{tikzpicture}
\caption{Smallest 100 eigenvalues $\lambda_{i}, i=1,\dots,100$ of the Maxwell eigenvalue problem computed with $M_h$ and $\widetilde{M}_h$. The spurious mode is marked with a red circle. The exact eigenvalues are not shown as they correspond to the ones computed with $\widetilde{M}_h$.}
 \label{fig:eigvals_spurious}
 \end{figure}

\begin{figure}
\centering
\includegraphics[width=0.3\textwidth]{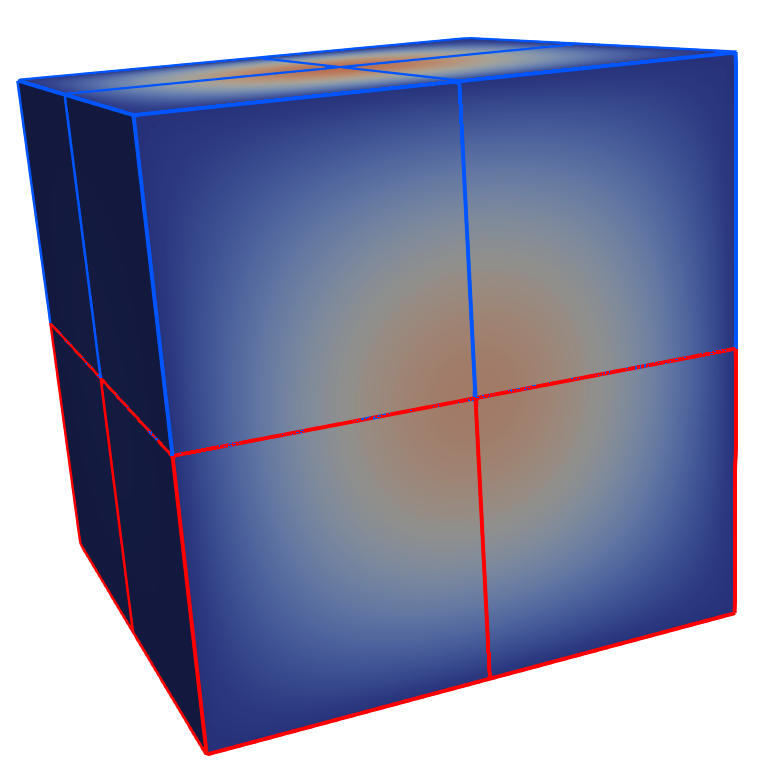}
\includegraphics[width=0.3\textwidth]{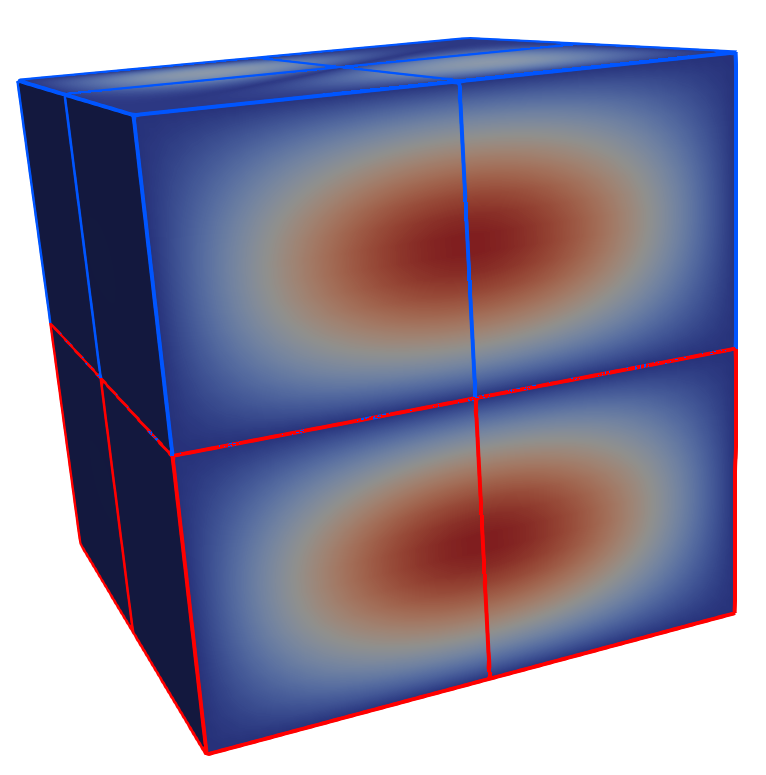}
\includegraphics[width=0.3\textwidth]{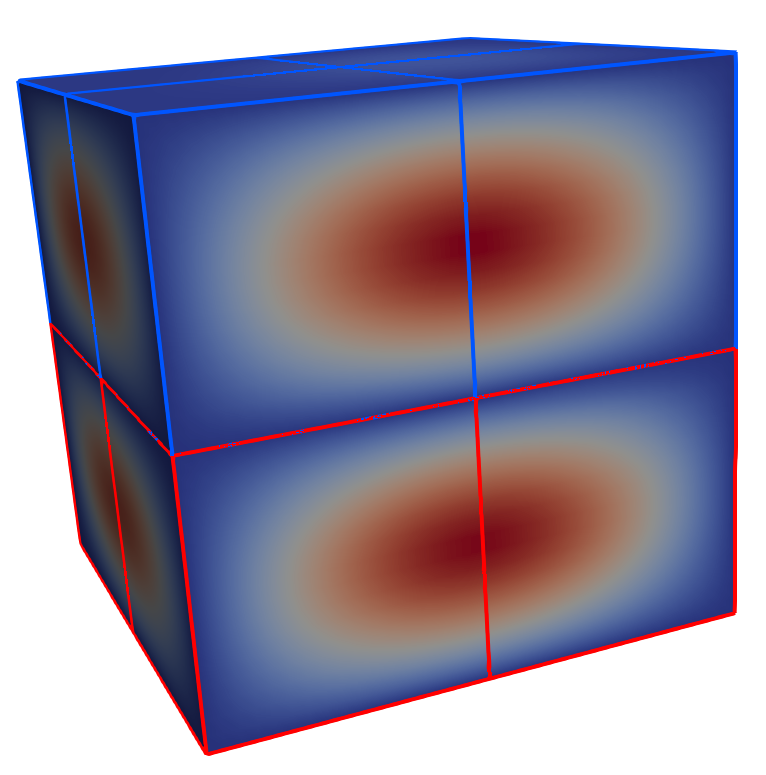}
\caption{Visualisation of some eigenmodes of the Maxwell eigenvalue problem.}
\label{fig:cube-modes}
\end{figure}

\begin{figure}
\centering
\includegraphics[width=0.3\textwidth]{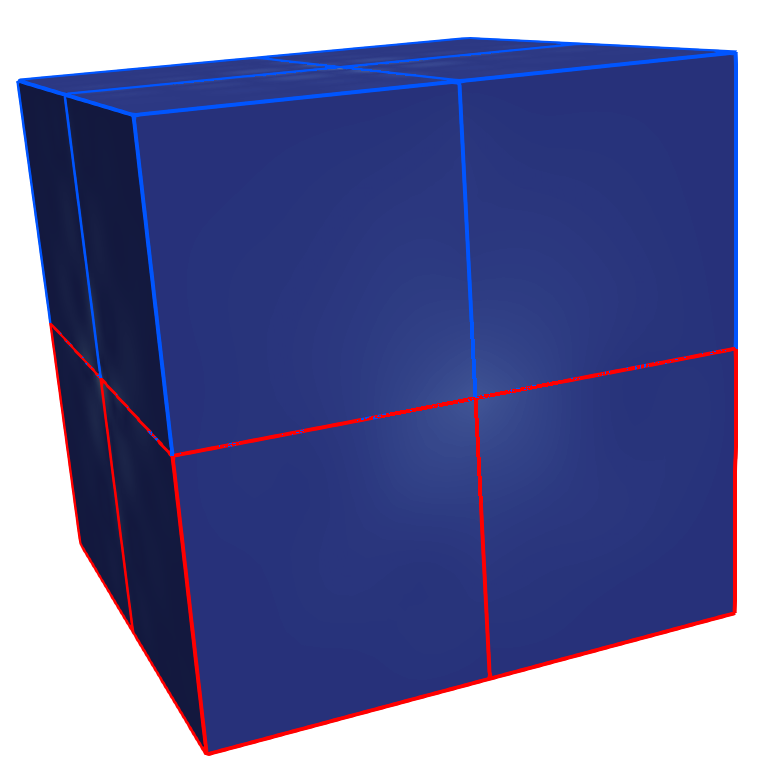}
\includegraphics[width=0.3\textwidth]{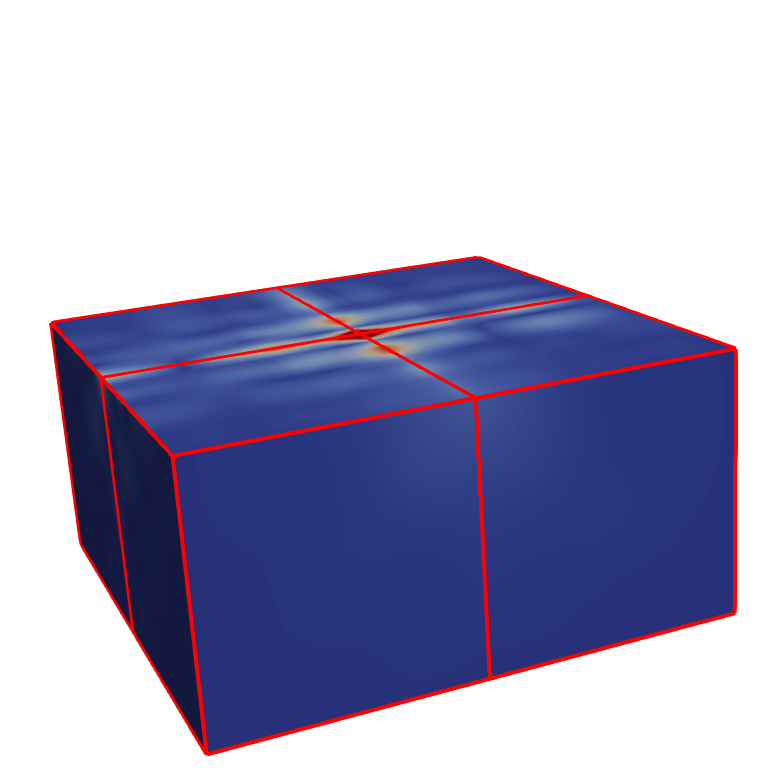}
\caption{Spurious mode introduced when using $M_h$ as space for the Lagrange multiplier.}
\label{fig:cube-spurious_mode}
\end{figure}

To analyse whether the discrete space $\widetilde{M}_h$ introduced in \cref{sec:mortar-gauged} is the appropriate choice in the case of mortaring we solve the Maxwell eigenvalue problem \cref{eq:eigvalprob} for a mortar problem in the domain \cref{fig:problemdomain_cube_conf} with Dirichlet boundary conditions. The resulting eigenvalues for a computation with degree $p=4$ in both $\Omega_1$ and $\Omega_2$ and $64$ elements
are shown in \cref{fig:eigvals_spurious}. A spurious mode is introduced when using the unmodified multiplier space $M_h$. When using $\widetilde{M}_h$ as multiplier space no spurious modes are introduced and the computed eigenvalues correspond to the exact ones. 
Some eigenmodes of the eigenvalue problem can be seen in \cref{fig:cube-modes}. When applying tree-cotree gauging without modifying the space $M_h$ we obtain a spurious mode for each internal vertex at the mortar interface. An example of such a spurious mode can be seen in \cref{fig:cube-spurious_mode}. When using the modified space $\widetilde{M}_h$ we get no spurious modes as expected.
\subsection{Magnetostatic Source Problem}
\begin{figure}
\centering
	\includegraphics[height=0.3\linewidth]{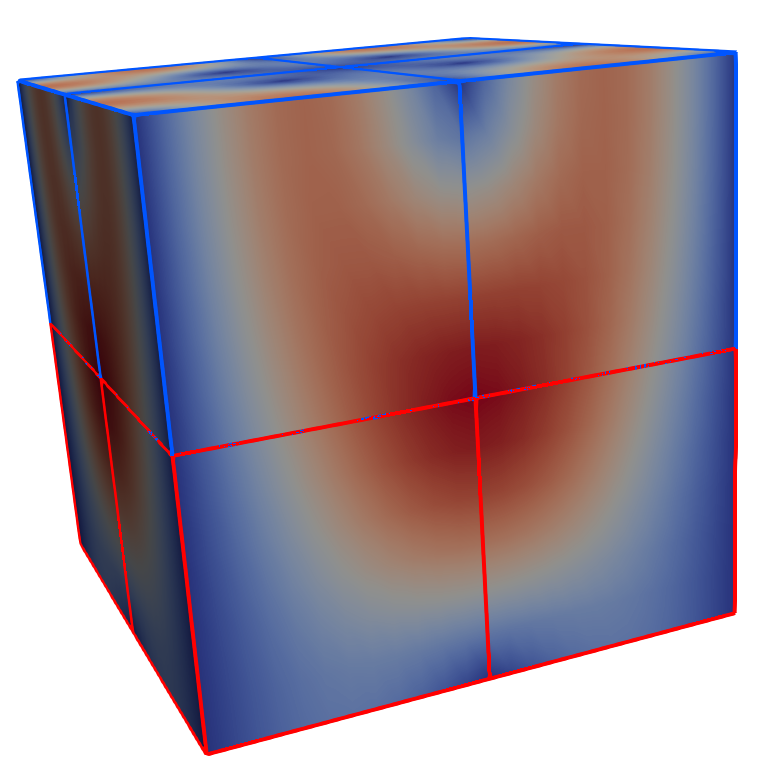}
	\includegraphics[height=0.3\linewidth]{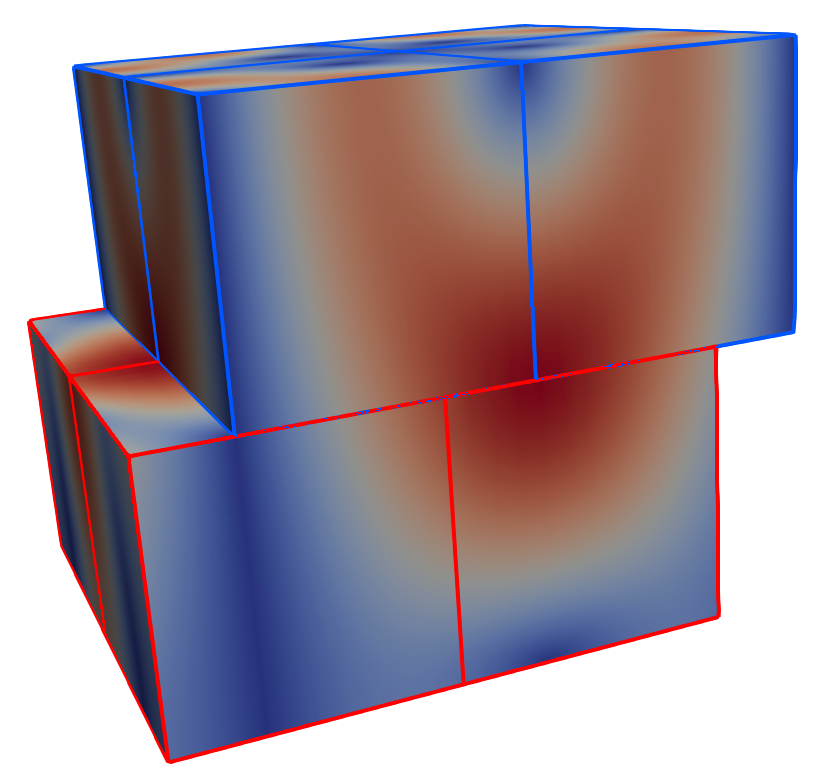}
\caption{Magnitude $|\vecs{B}|$ of the magnetic flux density of the test problem, normalised for the visualisazion.} \label{fig:Bfield_test}
\end{figure}

To analyse the convergence of the solution of the magnetostatic source problem using mortaring and tree-cotree gauging we solve \cref{eq:magpoteq,eq:neumann,eq:dirichlet} on the same domain shown in \cref{fig:problemdomain_cube} with periodic boundary conditions in $x$-direction and Dirichlet boundary conditions on all other boundaries. We use the manufactured solution
\begin{equation}
 \begin{aligned}
 \vecs{A}_{\mathrm{src}} \left( x,y,z\right)
 =
    \begin{bmatrix}
      \sin\left( y \right) \sin\left( \frac{z}{2} \right)  \\
      \sin\left( x \right) \sin\left( \frac{z}{2} \right)  \\
      \sin\left( x \right) \sin\left( y   \right) 
    \end{bmatrix}.
\end{aligned}
\end{equation}
This leads to the magnetic flux density
\begin{equation}
 \begin{aligned}
   \vecs{B}_{\mathrm{ana}} \left( x,y,z\right)
   =
    \begin{bmatrix}
      \sin\left( x \right) \cos\left( y \right) - \frac{1}{2}\left(\sin\left( x \right) \cos\left( \frac{z}{2} \right)\right)  \\
      - \cos\left( x \right) \sin\left( y \right) + \frac{1}{2}\left(\sin\left( y \right) \cos\left( \frac{z}{2} \right)\right)  \\
      \cos\left( x \right) \sin\left( \frac{z}{2} \right) - \cos\left( y \right) \sin\left( \frac{z}{2} \right)
    \end{bmatrix}. \label{eq:B}
 \end{aligned}
\end{equation}
which is induced by the electric current density and magnetisation
\begin{align}
 \vecs{J}_{\mathrm{src}} \left( x,y,z\right) &=
    \begin{bmatrix}
        \frac{5}{4} \sin\left( y \right) \sin\left( \frac{z}{2} \right) \\
        \frac{5}{4} \sin\left( x \right) \sin\left( \frac{z}{2} \right) \\
        2 \sin\left( x \right) \sin\left( y \right)
    \end{bmatrix},\\
 \vecs{M}\left( x,y,z\right) &= \vecs{0}.
\end{align}
The magnetic flux density \cref{eq:B} is shown in \cref{fig:Bfield_test}.
\Cref{fig:convergence_src} shows the convergence of the relative error 
\begin{equation}
 \epsilon = \frac{\sqrt{\int\limits_{\Omega}|\vecs{B}_{\mathrm{num}}-\vecs{B}_{\mathrm{ana}}|^{2} \operatorname{d}V}}{\sqrt{\int\limits_{\Omega}|\vecs{B}_{\mathrm{ana}}|^{2} \operatorname{d}V}}
\end{equation}
in the domain $\Omega = \overline{\Omega}_{1} \cup \overline{\Omega}_{2}$ of the numerically computed magnetic flux density $\vecs{B}_{\mathrm{num}}$ for the case of conforming and non-conforming patches. Here, mortaring is used with tree-cotree gauging using the multiplier space $\widetilde{M}_h$ introduced in \cref{sec:mortar-gauged}. We solve both for a domain with conforming patches, as in \cref{fig:problemdomain_cube_conf}, and for a domain with non-conforming patches obtained by shifting $\Omega_2$ by $\delta=1$ in $x$-direction, as in \cref{fig:problemdomain_cube_nonconf}. We observe that the solution converges with order $p$ both for the conforming and the non-conforming configuration.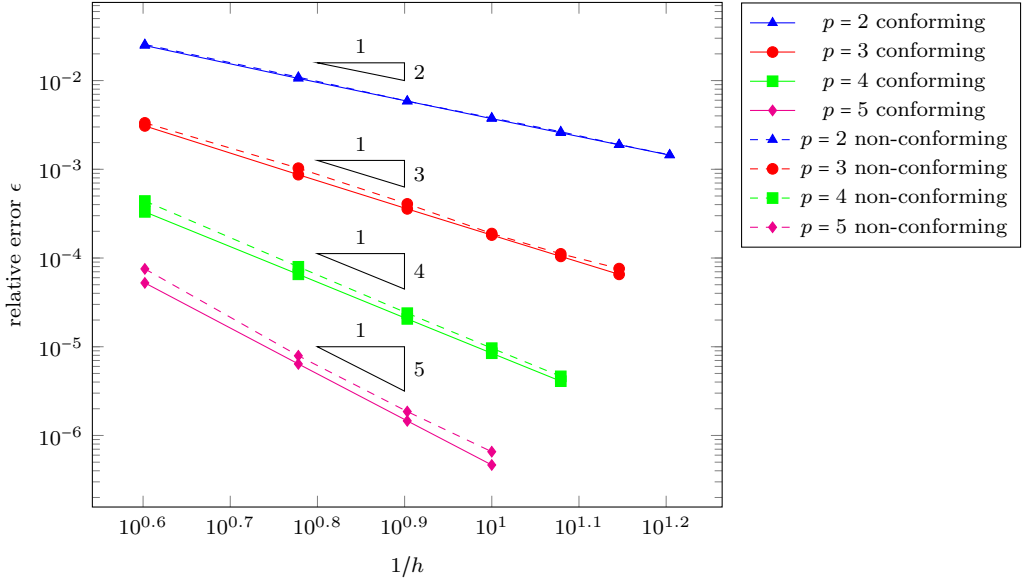
\begin{figure}
\centering
\begin{tikzpicture}[]
           \begin{loglogaxis}[width=0.6\columnwidth, height=0.5\columnwidth, ylabel={relative error $\epsilon$},xlabel={$1/h$}, ylabel near ticks,tick label style={font=\footnotesize},  legend style={font=\footnotesize, }, label style={font=\footnotesize},every x tick scale label/.style={at={(1,0)},anchor=north,yshift=-5pt,inner sep=0pt},legend pos=outer north east]
                \addplot [blue, mark=triangle*, mark options={solid}] table [x index=0, y index=5, col sep=comma] {images/convergence_cube_wholedom_rot0pi_nsub_ndof_relerr_deg_3_4_5.csv};
                \addplot [red, mark=*, mark options={solid}] table [x index=0, y index=6, col sep=comma] {images/convergence_cube_wholedom_rot0pi_nsub_ndof_relerr_deg_3_4_5.csv};
                \addplot [green, mark=square*, mark options={solid}] table [x index=0, y index=7, col sep=comma] {images/convergence_cube_wholedom_rot0pi_nsub_ndof_relerr_deg_3_4_5.csv};  
                \addplot [magenta, mark=diamond*, mark options={solid}] table [x index=0, y index=8, col sep=comma] {images/convergence_cube_wholedom_rot0pi_nsub_ndof_relerr_deg_3_4_5.csv};    
                \addplot [dashed, blue, mark=triangle*, mark options={solid}] table [x index=0, y index=5, col sep=comma] {images/convergence_cube_wholedom_rot1_nsub_ndof_relerr_deg_3_4_5.csv};
                \addplot [dashed, red, mark=*, mark options={solid}] table [x index=0, y index=6, col sep=comma] {images/convergence_cube_wholedom_rot1_nsub_ndof_relerr_deg_3_4_5.csv};
                \addplot [dashed, green, mark=square*, mark options={solid}] table [x index=0, y index=7, col sep=comma] {images/convergence_cube_wholedom_rot1_nsub_ndof_relerr_deg_3_4_5.csv};    
                \addplot [dashed, magenta, mark=diamond*, mark options={solid}] table [x index=0, y index=8, col sep=comma] {images/convergence_cube_wholedom_rot1_nsub_ndof_relerr_deg_3_4_5.csv};    
                \draw (axis cs:10^0.8,10^-1.8)--(axis cs:10^0.9,10^-2.0)--(axis cs:10^0.9,10^-1.8) node[pos=0.5,anchor=west] {\footnotesize$2$}--cyclenode[pos=0.5,anchor=south] {\footnotesize$1$};
                \draw (axis cs:10^0.8,10^-2.9)--(axis cs:10^0.9,10^-3.2)--(axis cs:10^0.9,10^-2.9) node[pos=0.5,anchor=west] {\footnotesize$3$}--cyclenode[pos=0.5,anchor=south] {\footnotesize$1$};
                \draw (axis cs:10^0.8,10^-3.95)--(axis cs:10^0.9,10^-4.35)--(axis cs:10^0.9,10^-3.95) node[pos=0.5,anchor=west] {\footnotesize$4$}--cyclenode[pos=0.5,anchor=south] {\footnotesize$1$};
                \draw (axis cs:10^0.8,10^-5)--(axis cs:10^0.9,10^-5.5)--(axis cs:10^0.9,10^-5) node[pos=0.5,anchor=west] {\footnotesize$5$}--cyclenode[pos=0.5,anchor=south] {\footnotesize$1$};
                \legend{$p=2$ conforming\\$p=3$ conforming\\$p=4$ conforming\\$p=5$ conforming\\$p=2$ non-conforming\\$p=3$ non-conforming\\$p=4$ non-conforming\\$p=5$ non-conforming\\
                }
            \end{loglogaxis}
\end{tikzpicture} 
\caption{Convergence of the relative error $\epsilon$ of the magnetic flux density.}
 \label{fig:convergence_src}
 \end{figure}
\subsection{Electric Machine}
The simulation and optimisation of electric machines with conventional finite elements is well established in academia and industry since many decades, see e.g. \cite{Salon_1995aa}. Commonly, the rotation is realised by a cut in the air gap and a domain decomposition into rotor and stator is applied. In the machine community, many methods have been proposed for handling the resulting non-matching discretisations at the interface, e.g.~locked step, sliding surface, moving band, harmonic coupling, and mortaring, see for example \cite{De-Gersem_2004ad,Davat_1985aa,Lai_1992aa,Bouillault_2001aa,Buffa_1999aa}. However, in the finite element context alls those methods suffer from mesh issues. Therefore, IGA was recently suggested for machine simulation, e.g. \cite{Bontinck_2018ac,Bhat_2018aa,Friedrich_2020aa}. NURBS-based geometry representation promises a robust numerical method since rotor and stator domains share an exact circular interface whose non-matching knot vectors can be overcome by mortaring as explained previously, see \cref{fig:images_mortar_interface}. In the following we demonstrate mortaring with tree-cotree gauging for the simulation of a six-pole permanent magnet synchronous machine. For the description of geometry and material coefficients see \cite[Chapter V.A]{Bontinck_2018af}.

\begin{figure}
\centering
 \begin{subfigure}{0.45\textwidth}
    \centering 
	\includegraphics[width=\textwidth]{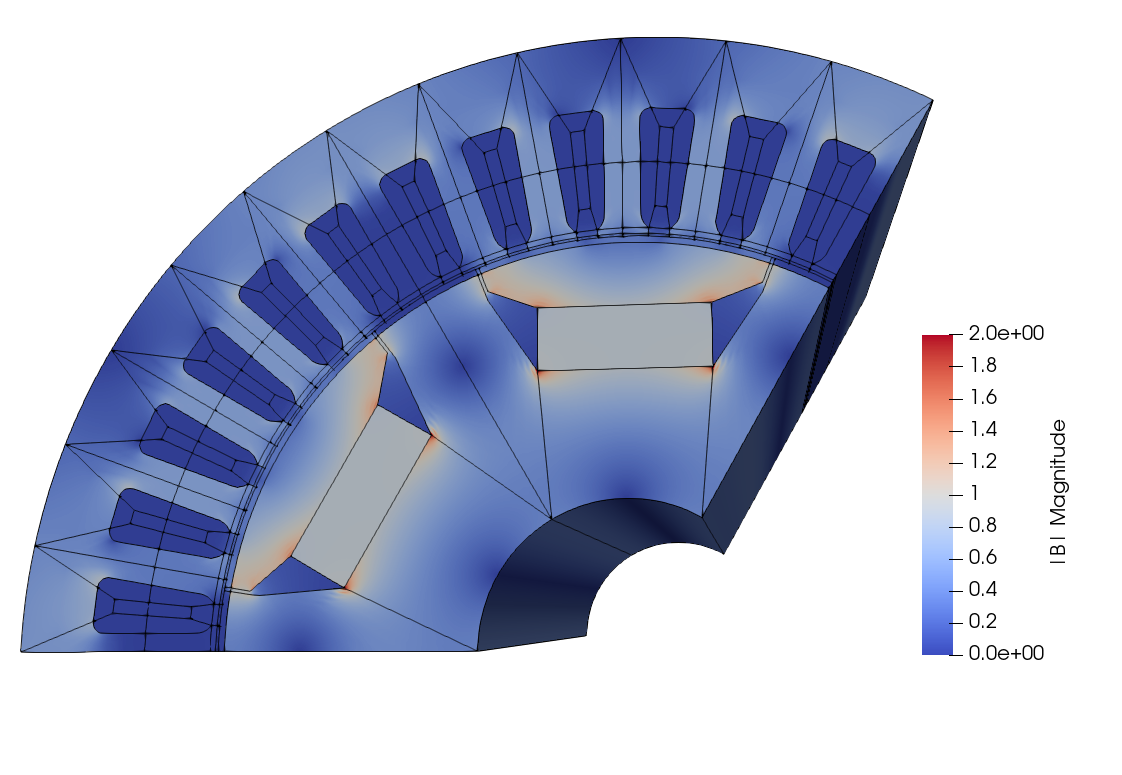}
 \caption{rotation angle $\alpha=0$}
\label{fig:machine_flux_rot0}
 \end{subfigure}
 \begin{subfigure}{0.45\textwidth}
    \centering 
	\includegraphics[width=\textwidth]{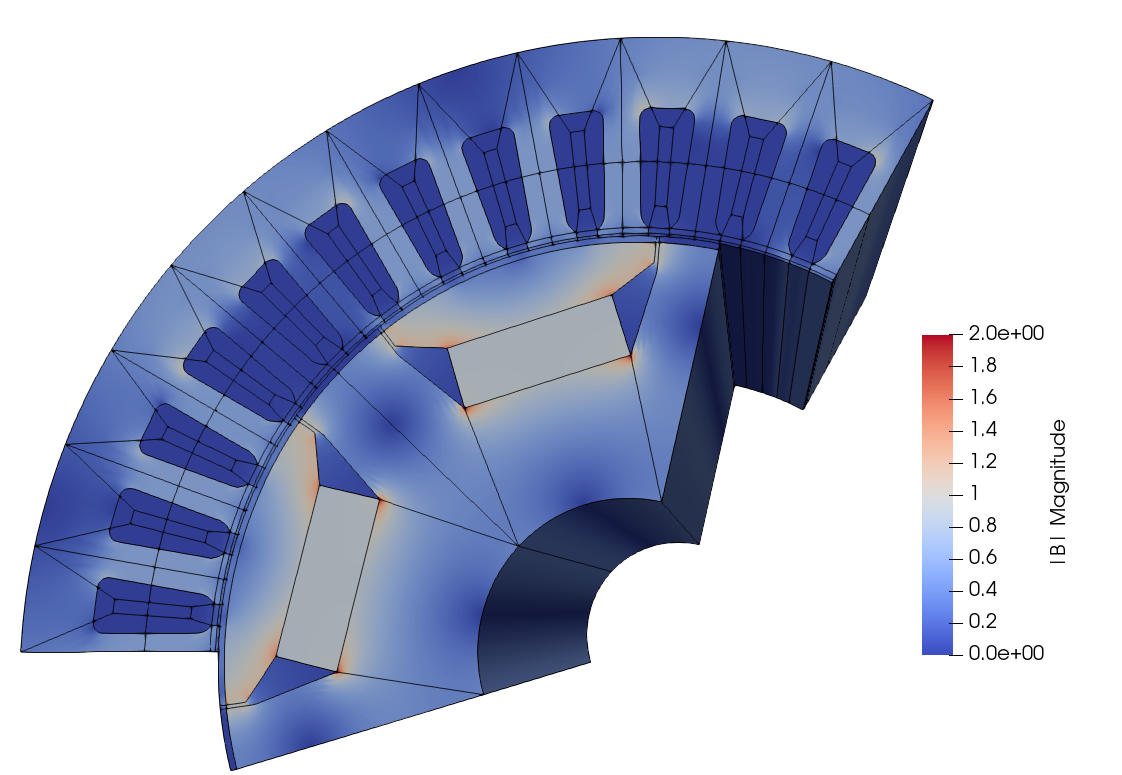}
 \caption{rotation angle $\alpha=\nicefrac{\pi}{11}$}
\label{fig:machine_flux_rot1}
 \end{subfigure}
\caption{Magnetic flux density in two poles of a six pole permanent magnet synchronous machine for different rotation angles.}
\label{fig:machine_flux}
\end{figure}

The magnetic flux density in one pole pair of the machine computed for two different rotor configurations is shown in \cref{fig:machine_flux}. The geometry is represented by $24$ rotor and $228$ stator patches. The total number of degrees of freedom of the linear equations system using tree-cotree gauging is $N_{\mathrm{dof}}=50330$. Only the coupling matrices , i.e. the $\mathbf{G}$ matrices in \cref{eq:matrix-gauged}, depend on the rotation angle of the machine, and thus, a rotating machine can be evaluated efficiently. The assembly of the coupling matrices takes about $\SI{5}{\percent}$ of the total assembly time and the creation of the tree and the cotree takes less than $\SI{2}{\percent}$ of the computational time.

\section{Conclusions}\label{sec:conclusions}
We have introduced the tree-cotree decomposition for curl-conforming spline spaces used in IGA. We have shown that using the correspondence of spline basis functions with edges of the control mesh, the existing tree-cotree decomposition algorithms from lowest order finite elements can be applied without changes, also for non-contractible domains. We have also extended the tree-cotree decomposition to multi-patch domains with non-conforming discretisations, glued together with the mortar method. The tree-cotree decomposition requires a simple modification in the dependent subdomain, to include all the degrees of freedom on the mortar interface in the cotree. We have found that, to couple mortaring and tree-cotree gauging while avoiding spurious solutions, it is also necessary to enrich the multiplier space with gradient functions of the internal vertices on interface. The numerical tests show that the enrichment of the multiplier space along with tree-cotree gauging provides a high-order method which is inf-sup stable and spurious-free. We have also presented an example of a rotating permanent magnet synchronous machine to show the potential of the method for the simulation of electric machines.
 
\appendix

\section{Analysis of the Discrete Kernel} \label{sec:appendix-kernel}
We analyse here the relation between discrete kernel $K_h$ and the spaces $\nabla X^0_{h,0}$ and $\nabla X^0_{h}$ in the setting of \cref{sec:mortar-gauged}. First of all, we can characterise $K_{h}$ analogously to the proof of Proposition~4.5 in \cite{Buffa_2020aa}: for any function $\mathbf{u}_h \in K_h$ it holds that $\nabla \times \mathbf{u}_h = \mathbf{0}$ in $\Omega_1$, and therefore there exists $\phi_h \in X_h^0$, but not necessarily in $X_{h,0}^0$, such that $\mathbf{u}_h = \nabla \phi_h$. Moreover, any function in $\nabla X^0_{h,0}$ is clearly in $K_h$, and therefore $\nabla X^0_{h,0} \subset K_{h} \subset \nabla X^0_{h}$. To prove that in general both inclusions are strict, we need several technical results.

Let $\mathbf{u}_h = \nabla \phi_h \in K_h$, with $\phi_h \in X_{h}^0$, and let us introduce the multi-patch trace space (see also \cite{Buffa_2019ac})
\[
W_h^0 = \gamma_{\Gamma_{\mathrm{int}}}^0 (V_h^0) = \{ \psi_h \in H^{1/2}_{00}(\Gamma_{\mathrm{int}}) : \psi_h |_{\Gamma^i} \in S_p^0(\Gamma^i; \partial \Gamma^i \cap \partial \Gamma_{\mathrm{int}}) \}.
\]
Then, since the trace operators and the differential operators commute, for $\mathbf{u}_h \in K_h$ it holds that
\[
{\gamma}^1_{\Gamma_{\mathrm{int}}} \mathbf{u}_h = 
{\gamma}^1_{\Gamma_{\mathrm{int}}} (\nabla \phi_h) = 
\nabla_\Gamma (\gamma^0_{\Gamma_{\mathrm{int}}} \phi_h) = \nabla_\Gamma \psi_h,
\]
for some discrete function $\psi_h \in W^0_h$. Let us also introduce the one-dimensional interfaces $\Lambda^{ij} = \partial \Gamma^i \cap \partial \Gamma^j$ for $1 \le i < j \le N_1$. By definition of $K_h$, and integrating by parts, we obtain
\begin{align}
(\nabla_\Gamma \psi_h, \boldsymbol{\mu}_h)_{\Gamma_{\mathrm{int}}} = 
- \sum_{i=1}^{N_1} (\psi_h, \nabla_\Gamma \cdot \boldsymbol{\mu}_h)_{\Gamma^i} + 
\sum_{i=1}^{N_1-1}\sum_{j=i+1}^{N_1} (\psi_h, \boldsymbol{\mu}_h |_{\Gamma^i} \cdot \mathbf{n}^i + \boldsymbol{\mu}_h |_{\Gamma^j} \cdot \mathbf{n}^j)_{\Lambda^{ij}} = 0 \quad \forall \boldsymbol{\mu}_h \in M_h. \label{eq:int-by-parts}
\end{align}
Analysing the expressions involving $\boldsymbol{\mu}_h$, from the definition of $M_h$ and the spaces introduced in \cref{sec:traces}, it holds that
\begin{align}
& \nabla_\Gamma \cdot \boldsymbol{\mu}_h |_{\Gamma^i} \in S^2_{p-1}(\Gamma^i), \text{ for } i = 1, \ldots, N_1, \label{eq:mu1}\\
& (\boldsymbol{\mu}_h \cdot \mathbf{n}^i)|_{\Lambda^{ij}} \in S^1_{p-1}(\Lambda^{ij}), \text{ for } 1\le i < j \le N_1, \label{eq:mu2}
\end{align}
Moreover, for each $\nu_h \in S_{p-1}^1(\Lambda^{ij})$ there exists a lifting $\boldsymbol{\eta}_h \in S_{p-1}^{1^*}(\Gamma^i)$ such that $\nabla_\Gamma \cdot \boldsymbol{\eta}_h = 0$ on $\Gamma^i$, and $\boldsymbol{\eta}_h \cdot \mathbf{n}^i |_{\Lambda^{ij}} = \nu_h$.

\begin{figure}
\centering
\includegraphics[width=0.4\textwidth,trim=2cm 1cm 1cm 1cm, clip]{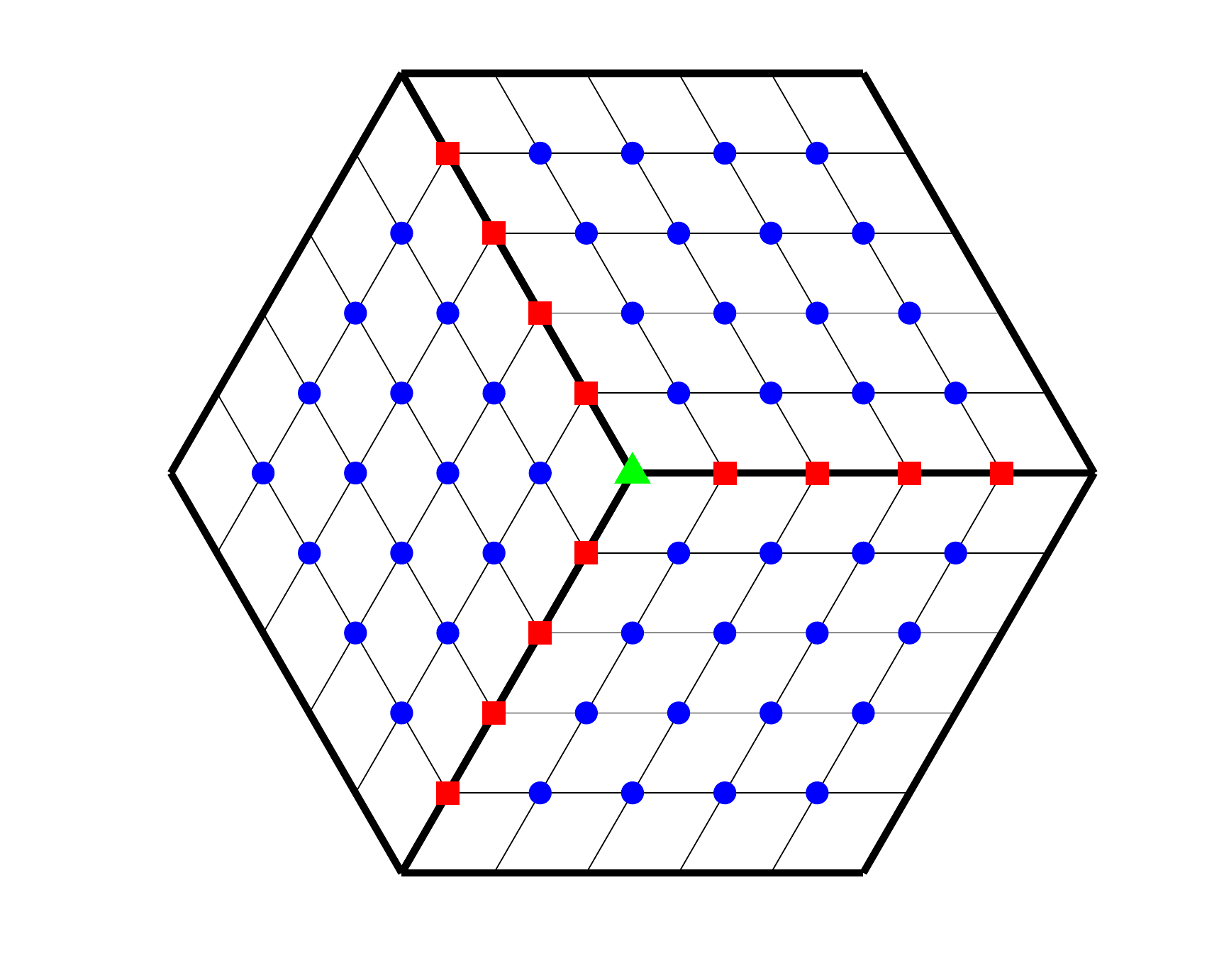}
\caption{Visualisation of the position of control points defining $\Gamma_{\mathrm{int}}$: internal to a patch (blue dots), on an edge (red squares) and on a vertex (green triangle).}
\label{fig:control_points_for_proof}
\end{figure}

As we are imposing boundary conditions strongly except on $\Gamma_{\mathrm{int}}$, the basis functions of $X^0_h$ can be associated to the control points of the control mesh, except those corresponding to $\partial \Omega_1 \Setminus \Gamma_{\mathrm{int}}$. Moreover, the internal control points correspond to functions in $X^0_{h,0}$, whose gradient is always in $K_h$. Therefore, it is enough to analyse the basis functions corresponding to control points defining $\Gamma_{\mathrm{int}}$, but not on $\partial \Gamma_{\mathrm{int}}$, and among these we have to distinguish three cases: control points internal to a patch of $\Gamma_{\mathrm{int}}$, points on an edge between two patches, and points on a vertex between three or more patches. They respectively correspond to blue circles, red squares and green triangles in \cref{fig:control_points_for_proof}.

For the first case (blue circles), and using the notation above for $\phi_h \in X^0_h$ and $\psi_h = \gamma^0_{\Gamma_{\mathrm{int}}} \phi_h$, it holds that $\mathrm{supp} \, \psi_h \subset \Gamma^i$ for some $i$, and $\psi_h|_{\Lambda^{ij}}$ = 0 for every $i,j$. Therefore, $\psi_h|_{\Gamma^i} \in S_p^0(\Gamma^i;\partial \Gamma^i)$. Let us assume that $\phi_h \in K_h$, which implies that \eqref{eq:int-by-parts} holds true. Using \eqref{eq:mu1}, the fact that the surface divergence operator is surjective from $S_{p-1}^{1^*}(\Gamma^i)$ into $S_{p-1}^2(\Gamma^i)$, and the inf-sup condition plus the equal dimension of $S_p^0(\Gamma^i;\partial \Gamma^i)$ and $S_{p-1}^2(\Gamma^i)$, we have that $\psi_h = 0$, which is a contradiction, and therefore $\phi_h \not \in K_h$.

Similarly, for the second case (red squares) we have $\mathrm{supp}\, \psi_h \subset \Gamma^i \cup \Gamma^j$ for some $i,j$, and $\psi_h |_{\Lambda^{ij}} = 0$. Let us assume that $\phi_h \in K_h$, and therefore that \eqref{eq:int-by-parts} holds. The divergence free liftings explained above yield $(\psi_h, \nu_h)_{\Lambda^{ij}} = 0$ for every $\nu_h \in S^1_{p-1}(\Lambda^{ij})$. But the inf-sup condition and equal dimension between the univariate spaces $S^0_p(\Lambda^{ij};\partial \Lambda^{ij})$ and $S_{p-1}^1(\Lambda^{ij})$ implies that $\psi_h|_{\Lambda^{ij}} = 0$, which is again a contradiction, and $\phi_h \not \in K_h$.

For the third case (green triangle) it is not possible in general to prove that the functions are not in $K_h$. As a consequence, the dimension of the kernel space is at most
\[
\dim X^0_{h,0} \le \dim K_h \le \dim X^0_{h,0} + \# Z_{\Gamma_{\mathrm{int}}}.
\]
The numerical tests presented in \cref{sec:mortar-gauged} lead us to conjecture that the second inequality is in fact an equality, but we have not been able to prove that the equality is always true.

\section*{Acknowledgement}
This work is supported by the Graduate School CE within the Centre for Computational Engineering at Technische Universität Darmstadt, by the Swiss National Science Foundation via the project HOGAEMS n.200021\_188589, by the German Research Foundation via the project SCHO 1562/6-1, and by the Defense Advanced Research Projects Agency (DARPA), under contract HR0011-17-2-0028. The views, opinions and/or findings expressed are those of the author and should not be interpreted as representing the official views or policies of the Department of Defense or the U.S. Government.

\end{document}